\documentclass[]{amsart}
\counterwithin{equation}{section}
\usepackage[scale=0.75]{geometry}

\usepackage{bbm}
\usepackage{amssymb}
\usepackage{amsmath}
\usepackage[all]{xy}
\usepackage{tikz}
\usepackage{amsfonts}
\usepackage{mathrsfs}
\usepackage{latexsym}
\usepackage{graphicx}
\usepackage{enumerate}
\usepackage{amscd,amssymb,amsmath,amsbsy,amsthm}
\definecolor{linkred}{rgb}{0.7,0.2,0.2}
\definecolor{linkblue}{rgb}{0,0.2,0.6}
\definecolor{linkgreen}{rgb}{0,0.6,0.2}
\usepackage[colorlinks,plainpages,backref,
    linkcolor=linkblue,
    citecolor=linkgreen,
    urlcolor=linkred]{hyperref}

\newtheorem{Th}{Theorem}[section]
\newtheorem{Lemma}[Th]{Lemma}
\newtheorem{Def}[Th]{Definition}
\newtheorem{Coro}[Th]{Corollary}
\newtheorem{Prop}[Th]{Proposition}
\newtheorem{Eg}[Th]{Example}
\newtheorem{Rmk}[Th]{Remark}

\allowdisplaybreaks

\newtheorem{Thm}{Theorem}[section]

\def\four#1#2#3#4{\frac{\theta(#1)\theta(#2)}{\theta(#3)\theta(#4)}}
\def\P{\mathsf{P}}
\def\Q{\mathsf{Q}}
\def\calQ{\mathcal{Q}}

\def\bfE{\mathbf{E}}

\def\v{{\mathsf{v}}}
\def\la{\lambda}
\def\de{\delta}
\def\al{\alpha}
\def\alv{{\alpha^\v}}
\def\be{\beta}

\def\-{{-}}
\def\+{\texttt{+}}
\def\id{{\mathsf{id}}}
\def\pt{{\mathsf{pt}}}
\def\dyn{{\mathsf{d}}}
\def\frakE{{\mathfrak{E}}}
\def\cancel#1{{\not{\color{gray}\!#1}}}

\newcommand{\BPD}[2][1pc]{%
\setlength{\unitlength}{#1}
\def\FF{%
    \qbezier(0.5,0)(0.5,0.2)(0.5,0.2)
    \qbezier(1,0.5)(0.8,0.5)(0.8,0.5)
    \qbezier(0.8,0.5)(0.5,0.5)(0.5,0.2)
}
\def\JJ{%
    \qbezier(0.5,1)(0.5,0.8)(0.5,0.8)
    \qbezier(0,0.5)(0.2,0.5)(0.2,0.5)
    \qbezier(0.5,0.8)(0.5,0.5)(0.2,0.5)
    }
\def\II{%
    \qbezier(0.5,0)(0.5,0.5)(0.5,1)}%
\def\HH{%
    \qbezier(0,0.5)(0.5,0.5)(1,0.5)}%
\def\XX{\II\HH}
\def\NN{%
    \qbezier(0.5,0)(0.5,0.3)(0.5,0.3)
    \qbezier(0.5,1)(0.5,0.7)(0.5,0.7)}%
\def\BPDfr##1##2{%
\begin{picture}(1,1)%
    \linethickness{0.08\unitlength}
    ##1
    \put(0.5,0.2){\makebox[0pc]{\makebox[\unitlength][c]{\color{black}\(##2\)}}}
    \thinlines
    \color{lightgray}%
    \put(0,0){\line(0,1){1}}%
    \put(1,0){\line(0,1){1}}%
    \put(0,0){\line(1,0){1}}%
    \put(0,1){\line(1,0){1}}%
\end{picture}}
\providecommand{\Op}[1][]{\BPDfr{}{##1}}
\let\O\Op
\providecommand{\Xp}[1][]{\BPDfr{\XX}{##1}}
\let\X\Xp
\providecommand{\Fp}[1][]{\BPDfr{\FF}{##1}}
\let\F\Fp
\providecommand{\Jp}[1][]{\BPDfr{\JJ}{##1}}
\let\J\Jp
\providecommand{\Hp}[1][]{\BPDfr{\HH}{##1}}
\let\H\Hp
\providecommand{\Ip}[1][]{\BPDfr{\II}{##1}}
\let\I\Ip
\providecommand{\Bp}[1][]{\BPDfr{\FF\JJ}{##1}}
\let\B\Bp
\def\+{\begin{picture}(1,1)
    \color{gray}
    \put(0.5,0){\line(0,1){1}}
    \put(0,0.5){\line(1,0){1}}
    \put(0.5,0.5){\circle*{0.2}}
\end{picture}}%
\def\+{\BPDfr{}{\color{gray}?}}
\def\M##1{\begin{picture}(1,1)%
    \put(0,0.2){\makebox[\unitlength]{\(##1\)}}
\end{picture}}%
\begin{array}{@{\,}c@{\,}}
{\def\arraystretch{0}
\setlength{\arraycolsep}{0pc}
\color{teal}
\begin{array}{@{}l@{}}%
#2\end{array}}
\end{array}}

\def\redstrings{%
    \def\b{--++(.25,0)}
    \def\u{--++(.25,0)--++(1,1)--++(.25,0)}
    \def\d{--++(.25,0)--++(1,-1)--++(.25,0)}
    \def\h{--++(1.5,0)}
    \def\e{--++(.25,0)}}
\def\bluestrings{%
    \def\b{--++(.25,0)}
    \def\u{--++(.20,0)--++(1,1)--++(.30,0)}
    \def\d{--++(.30,0)--++(1,-1)--++(.20,0)}
    \def\h{--++(1.5,0)}
    \def\e{--++(.25,0)}}

\newcommand{\Rui}[1]{\textcolor{green}{$[$ Rui: #1 $]$}}

\begin{document}

\title{Combinatorial Aspects of Elliptic Schubert Calculus}

\author{Cristian Lenart}
\address{Department of Mathematics and Statistics,
State University of New York at Albany}
\email{clenart@albany.edu}

\author{Rui Xiong}
\address{Department of Mathematics and Statistics, University of Ottawa, 150 Louis-Pasteur, Ottawa, ON, K1N 6N5, Canada}
\email{rxion043@uottawa.ca}

\author{Changlong Zhong}
\address{Department of Mathematics and Statistics,
State University of New York at Albany}
\email{czhong@albany.edu}


\begin{abstract}
The main goal of this paper is to extend two fundamental combinatorial results in Schubert calculus on flag manifolds from equivariant cohomology and $K$-theory to equivariant elliptic cohomology. The foundations of elliptic Schubert calculus were laid in a few relatively recent papers by Rim\'anyi, Weber, and Kumar. They include the recursive construction of elliptic Schubert classes via generalizations of the cohomology and $K$-theory push-pull operators and the study of the corresponding Demazure algebra. We derive a {Billey-type formula} for the localization of elliptic Schubert classes (for partial flag manifolds of arbitrary type) and a {pipe dream model} for their polynomial representatives in the case of type $A$ flag manifolds. The latter extends the pipe dream model for double Schubert and Grothendieck polynomials. We also study the degeneration of elliptic Schubert classes to $K$-theory, which recovers the corresponding classical formulas.  
\end{abstract}

\maketitle


\section{Introduction}\label{intro}

There are three kinds of one-dimensional algebraic groups: the additive group $\mathbb{G}_a$, the multiplicative group $\mathbb{G}_m$, and elliptic curves. The corresponding generalized cohomology theories are cohomology, $K$-theory, and elliptic cohomology. Throughout this paper, we work with (torus) equivariant elliptic cohomology.

From the construction, typical elements in elliptic cohomology are elliptic functions in line bundles, which makes elliptic cohomology less understood than cohomology and $K$-theory. 
By the seminal work of Aganagic and Okounkov \cite{AO21}, elliptic cohomology provides a framework to control the generating functions in the $K$-theoretic curve countings. It gains importance due to its close connection to an important duality called ``3D mirror symmetry'' from mathematical physics; see \cite{Smi24} for an introduction. 

Cohomological and $K$-theoretic Schubert calculus, including their equivariant versions, have been extensively studied from various perspectives for many years. In contrast, elliptic Schubert calculus is a relatively recent subject. Its foundations were laid in~\cite{RW1,KRW,LZZ,Rim21}, and they include the recursive construction of elliptic Schubert classes via generalizations of the cohomology and $K$-theory push-pull operators and the study of the corresponding Demazure algebra. However, elliptic Schubert calculus remains underdeveloped, particularly in terms of its combinatorial aspects. 
In this paper, we derive the analogues in equivariant elliptic cohomology of two fundamental combinatorial formulas in Schubert calculus on partial flag manifolds: the \emph{Billey-type formula} for the localization of Schubert classes, and the \emph{pipe dream model} for their polynomial representatives in type $A$. 

As noticed by Rim\'anyi and Weber \cite{RW1}, the failure of well-definedness of fundamental classes leads to a remarkable new feature of elliptic Schubert calculus: the dependence of \emph{dynamical parameters}. 
This feature tremendously increases the complexity, but at the same time unifies the corresponding formula in a more systematic way. 
For example, the parabolic Schubert classes can be obtained directly from Schubert classes on full flag manifolds by specializing dynamical parameters, which cannot be seen in $K$-theoretic or cohomological Schubert calculus. 

Let us turn to the statement of our theorems. 
In the present paper, we consider the following two functions: 
$$\P(x,y)=\four{x-y}{\hbar}{y+\hbar}{x},\qquad 
\Q(x,y)=\four{x+\hbar}{y}{y+\hbar}{x},$$
where $\theta$ is the Jacobi theta function and $\hbar$ is a formal variable. 
Let $G$ be a reductive group. 
We use $z_\al$ (resp. $\la_\alv$) to represent the equivariant parameter for a root $\al$ (resp., the dynamical parameter for a coroot $\alv$) of $G$.
We now state our elliptic Billey-type formula.

\begin{Thm}[Theorem \ref{th:Billey}]
Let $u,w\in W$. Let $u=s_{i_1}\cdots s_{i_\ell}$ be a reduced decomposition. 
The localization of the elliptic Schubert class $\bfE_w$ indexed by $w$ admits the following combinatorial formula:
$$\bfE_w(u)=
\sum_{J}\prod_{j=1}^\ell\begin{cases}
\Q(\la_{\check{\gamma}_{j}^{J}},z_{\beta_j}),& j\in J,\\
\P(\la_{\check{\gamma}_{j}^{J}},z_{\beta_j}),& j\notin J,
\end{cases}\qquad
\text{where }
\left\{\begin{aligned}
\epsilon_j&=\delta^{\mathsf{Kr}}_{j\in J},\\
\beta_j & = s_{i_1}\cdots s_{i_{j-1}}\alpha_{i_j},\\
\check{\gamma}_{j}^J& = 
s^{\epsilon_{\ell}}_{i_\ell}\cdots s^{\epsilon_{j+1}}_{i_{j+1}}\alpha^\v_{i_j}. 
\end{aligned}\right.
$$
with the sum over subwords $J\subset \{1,\ldots,n\}$ such that $w=s_{i_1}^{\epsilon_1}\cdots s_{i_\ell}^{\epsilon_\ell}$. 
\end{Thm}

The above formula generalizes the original one for cohomology Schubert classes given in~\cite{AJS, billey} and its $K$-theory analogue (for structure sheaves of Schubert varieties and their duals) in~\cite{Gra, Wil}, while it extends the Billey-type formula in hyperbolic cohomology derived in~\cite{LZ}. In fact, by taking the limit $q\to 0$, our formula recovers the formula for Segre motivic Chern (SMC) classes of Schubert cells~\cite{SZZ}, which implies the mentioned formula for structure sheaves of Schubert varieties by specialization (see, e.g., \cite{AMSS}). In this sense, our formula also generalizes the Billey-type formula for the Segre-Schwartz-MacPherson (SSM) classes of Schubert cells in cohomology~\cite{Su}.

In contrast to  Billey's original paper \cite{billey}, we do not use the Yang-Baxter equation in our proof; in fact, we do the opposite: we derive the Yang-Baxter property (of the corresponding dynamical $R$-matrices) from our Billey-type formula (Corollary \ref{coro:YBequation}). 
In addition, we derive the parabolic version (Theorem \ref{th:Billeypara}) of our formula by specializing dynamical parameters, after which it suffices to sum over a subset of subwords. 
If we restrict to type $A$, our Billey-type formula admits a wiring diagram presentation (Theorem \ref{th:BilleytypeA}). 
Using that, we give a combinatorial proof of an identity of 3D mirror symmetry by constructing a simple sign-reversing involution. 

Our second theorem is a combinatorial formula for polynomial representatives of equivariant elliptic Schubert classes (in the sense specified in Section~\ref{subs:polrep}) in type $A$. This formula extends the pipe dream model for double Schubert and Grothendieck polynomials.  We use the \emph{generic pipe dreams} of Knutson and Zinn-Justin~\cite{KZ}, which they use to compute the SSM and SMC classes of Schubert cells of full flag varieties in type $A$. Thus, our result also extends theirs to the elliptic case. 
A {generic pipe dream} is a tiling of the following $n\times n$ square grid by the following tiles:
$$
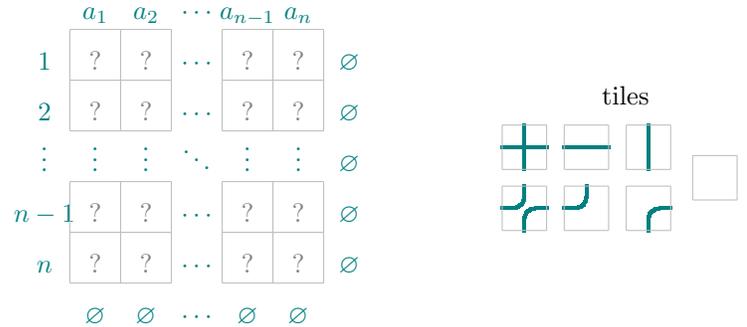
\BPD[1.6pc]{
\M{}\M{a_1}\M{a_2}\M{\cdots}\M{a_{n-1}}\M{a_n}\\
\M{1}\+\+\M{\cdots}\+\+\M{\varnothing}\\
\M{2}\+\+\M{\cdots}\+\+\M{\varnothing}\\
\M{\vdots}\M{\vdots}\M{\vdots}\M{\ddots}\M{\vdots}\M{\vdots}\M{\varnothing}\\
\M{n-1}\+\+\M{\cdots}\+\+\M{\varnothing}\\
\M{n}\+\+\M{\cdots}\+\+\M{\varnothing}\\
\M{}\M{\varnothing}\M{\varnothing}\M{\cdots}\M{\varnothing}\M{\varnothing}}\qquad\qquad 
\begin{array}{c}
\text{tiles}\\[1ex]
\begin{array}{@{}c@{}}
\BPD[1.4pc]{\X}\\[2ex]
\BPD[1.4pc]{\B}
\end{array}\,\,
\begin{array}{@{}c@{}}
\BPD[1.4pc]{\H}\\[2ex]
\BPD[1.4pc]{\J}
\end{array}\,\,
\begin{array}{@{}c@{}}
\BPD[1.4pc]{\I}\\[2ex]
\BPD[1.4pc]{\F}
\end{array}\,\,\,
\BPD[1.4pc]{\O}\,.
\end{array}
$$
For $w\in S_n$, we denote by $\mathsf{GPD}(w)$ the set of generic pipe dreams with $w(a_i)=i$, i.e., the $i$-th pipe from the left boundary goes to the $w(i)$-th position of the upper boundary. 

To get a polynomial representative, we need to associate with each tile an elliptic weight. 
As a feature of elliptic cohomology, each pipe is associated to an independent dynamical variable $\la_i$. The elliptic weight is determined by: the position of the corresponding tile, the indices of pipes inside, and another statistics called level (Definition \ref{def:level}). 
More precisely, we define: 
$$
\begin{array}{cc}
\BPD[1.2pc]{\X} & \Q(\la_a\-\la_b,y_j\-x_i) \\[2ex]
\BPD[1.2pc]{\B} & \P(\la_a\-\la_b,y_j\-x_i)
\end{array}\quad
\begin{array}{cc}
\BPD[1.2pc]{\H} & \Q(\la_a\-c\hbar,y_i\-x_i)
\\[2ex]
\BPD[1.2pc]{\J} & \P(\la_a\-c\hbar,y_i\-x_i)
\end{array}\quad
\begin{array}{cc}
\BPD[1.2pc]{\I} & \Q(c\hbar\-\la_b,y_i\-x_i)
\\[2ex]
\BPD[1.2pc]{\F} & \P(c\hbar\-\la_b,y_i\-x_i)
\end{array}
\quad 
\begin{array}{cc}
\BPD[1.2pc]{\O} & 1 
\end{array}$$
where $a$ (resp., $b$) is the index for the pipe from the lower (resp., left) boundary and $c$ is the level.

\begin{Thm}[Theorem \ref{thm:pipedreamofEw}]
For any $w\in S_n$, the weighted sum 
\begin{equation}\label{eq:thm:pipedreamofEw}
\mathcal{E}_w = \sum_{\pi\in \mathsf{GPD}(w)}\operatorname{weight}(\pi)
\end{equation}
is a polynomial representative of the elliptic Schubert class $\bfE_w$. 
\end{Thm}

The idea of the proof is to derive the pipe dream formula directly from the (type $A$ version) of our parabolic Billey-type formula. This idea is more conceptual and simpler than the original proofs of the pipe dream formulas for double Schubert and Grothendieck polynomials in~\cite{BJS, FS, FK}. 
More precisely, we realize the right-hand side of \eqref{eq:thm:pipedreamofEw} as a parabolic Billey-type formula in $S_{2n}$, and we use the fact that the corresponding localizations of elliptic Schubert classes satisfy the required recursion formula for our polynomial representatives. In the case of type $A$ Schubert classes in equivariant cohomology, the possibility of identifying their localizations with double Schubert polynomials (via pipe dreams) first appeared  in \cite{BR04}. 

The paper is organized as follows. 
In Section \ref{sec:ellSchcls}, we give the definition of elliptic Schubert classes and study their recursion.
Section \ref{sec:ellBilfor} is devoted to the Billey-type formula.
In Section \ref{sec:diatypeA}, we translate the formulas in type $A$ in terms of sub-wiring diagrams.
In Section \ref{ref:pipedream}, we establish the pipe dream model for a polynomial representative. 
Appendix \ref{sec:Kthylimit} explains the limit to $K$-theory.
In Appendix \ref{sec:comparing}, we explain the relation between our classes and those in \cite{RW1} and \cite{KRW}. 

\subsection*{Acknowledgment} C.L. was partially supported by the NSF grants DMS-1855592 and DMS-2401755, and C.Z. was partially supported by Simons Foundation Travel Support for Mathematicians TSM-00013828. 
C.L. and C.Z. would like to thank Gufang Zhao for the related collaboration, and Anders Buch, Allen Knutson and Richard Rim\'anyi for helpful discussions. 

\section{Preliminaries}

\subsection{Elliptic functions}
We recall the classical \emph{Jacobi theta function}
\begin{equation}\label{eq:thetadef}
\theta(u) = (x^{1/2}-x^{-1/2})\prod_{n>0}(1-q^nx)(1-q^n/x)\text{ where }
x=e^{2\pi i u}.
\end{equation}
The theta function is holomorphic in $u$ if  $|q|<1$, and it gives a meromorphic section of a certain line bundle over the elliptic curve $\mathbb{C}^\times/q^{\mathbb{Z}}$.  
In the present paper, we treat $\theta$ formally. 
More precisely, given a lattice $\Lambda$ and an element $u$, we could view 
\begin{equation}\label{eq:convtheta}
\theta(u)\in \mathbb{Q}[\Lambda]
[[q]], 
\mbox{\text{ where 
$\mathbb{Q}[\Lambda]=\bigoplus_{\lambda\in \Lambda} 
\mathbb{Q}\cdot e^{\pi i\lambda}$ is the group algebra}}.
\end{equation}
We  fix the following two functions: 
\begin{equation}\label{eq:defofPQ}
\P(x,y)=\four{x-y}{\hbar}{y+\hbar}{x},\qquad 
\Q(x,y)=\four{x+\hbar}{y}{y+\hbar}{x}.
\end{equation}
They  appear in the coefficients of elliptic Demazure-Lusztig operators and also in our combinatorial formulas. We record two identities.

\begin{Prop}
We have the following identities:
\begin{enumerate}[\quad \rm(1)]
    \item $\Q(x,y)\Q(y,x)=1$
    and $\P(x,y)+\Q(x,y)\P(y,x)=0$; 
    \item $\P(-\hbar,y)=1$ and $\Q(-\hbar,y)=0$;
    \item $\P(x,0)=1$ and $\Q(x,0)=0$. 
\end{enumerate}
\end{Prop}
The three sets of identities can be rewritten in terms of matrix identities:
\begin{align}
\label{eq:PQinv}
\left[\begin{matrix}
1 & 0 \\
\P(x,y) & \Q(x,y)
\end{matrix}\right]\cdot 
\left[\begin{matrix}
1 & 0 \\
\P(y,x) & \Q(y,x)
\end{matrix}\right]
& =\left[\begin{matrix}1&0\\0&1\end{matrix}\right];\\
\label{eq:PQnorm}
\left[\begin{matrix}
1 & 0 \\
\P(-\hbar,y) & \Q(-\hbar,y)
\end{matrix}\right]
=
\left[\begin{matrix}
1 & 0 \\
\P(x,0) & \Q(x,0)
\end{matrix}\right]
& =
\left[\begin{matrix}1&0\\1&0\end{matrix}\right].
\end{align}

\subsection{Root systems} Let $G$ be a reductive group with Borel subgroup $B$. Let $T$ be the maximal torus and $X^*(T), X_*(T)$ the set of characters and cocharacters, respectively. 
Let $\Sigma$ be  the set of simple roots, $\Phi$ the set of roots, and  $\Phi^+$ the set of positive roots. Let $W$ be the Weyl group. Let $P$ be a parabolic subgroup of $G$, which corresponds to a subset $\Sigma_P$ of $\Sigma$. Let $\Phi_P\subset \Phi$ be the subset of roots corresponding to $P$. Let $W_P$ be the subgroup generated by $s_\al, \al\in \Sigma_P$, and let $W^P$ be the set of minimal length representatives of left cosets of $W_P$ in $W$. 

\section{Elliptic Schubert Classes}
\label{sec:ellSchcls}

In this section, we review  the elliptic Schubert classes via the elliptic Demazure-Lusztig operators and compare them with the elliptic Schubert classes in \cite{RW1}.

\subsection{Twisted group algebras}
Let $T$ be the maximal torus of $G$, and let $\hbar$ be a formal parameter. We consider 
\begin{equation}\label{eq:defofQcal}
\mathcal{Q}=
\operatorname{Frac}\mathbb{Q}[\Lambda][[q]]
\text{ where }
\Lambda = X^*(T)\oplus X_*(T)\oplus
\mathbb{Z}\hbar. 
\end{equation}
For a character $\alpha\in X^*(T)$ (resp., a cocharacter $\beta^\v\in X_*(T)$), we denote by $z_\alpha$ (resp. $\la_{\beta^\v}$) the corresponding element in $\Lambda$. 
We think of $z$-variables as ``equivariant parameters'' and $\la$-variables as ``dynamical parameters.''  
The algebra $\mathbb{Q}[\Lambda]$ admits two commutative $W$-actions (and so does $\mathcal{Q}$): by acting on the $z$-variables and $\la$-variables, respectively; the latter is called the dynamical action. These actions are given by 
${}^{wv^\dyn}z_\alpha = z_{w\alpha} \text{ and }
{}^{wv^\dyn}\la_{\beta^\v} = \la_{v\beta^\v}$. 
Define the \emph{twisted group algebra}
\begin{equation}\label{eq:defofQcalW2}
\mathcal{Q}_{W^2} = \mathcal{Q}\rtimes \mathbb{Q}[W\times W^\dyn]. 
\end{equation}
It is a free left $\mathcal{Q}$-module with basis $\{\de_w\de_v^{\dyn}: w,v\in W\}$,  
with multiplication given by 
$$
a\,\de_w\de_v^{\dyn} \cdot 
a'\,\de_{w'}\de_{v'}^{\dyn}
= a\cdot {}^{wv^\dyn}a'\,\de_{ww'}\de^\dyn_{vv'}, \quad a,a'\in \calQ. $$
We can view it as the algebra generated by two $W$-actions on $\mathcal{Q}$ and 
the operators of multiplication by elements of $\mathcal{Q}$. 
Next, consider the algebra
\begin{equation}
\label{eq:defofQcalWstar}
\mathcal{Q}_W^* = \operatorname{Map}(W,\mathcal{Q})
\end{equation}
given by componentwise addition and multiplication. 
We denote by $f_w\in \mathcal{Q}_W^*$ the element defined by $f_w(v)=\delta_{w,v}^{\mathsf{Kr}}$.
In other words, for any $f\in \mathcal{Q}_W^*$, we can write 
$f=\sum_{w\in W} f(w)f_w \in \mathcal{Q}_W^*$. 

Let us briefly explain the geometric meaning of $\mathcal{Q}_W^*$.  
In this paper, we use the $K$-theory model of elliptic cohomology, i.e. we take the \emph{equivariant elliptic cohomology} to be
\begin{equation}\label{eq:defEGB}
\mathbf{E}ll_T(G/B)=K_T(G/B)[[q]].
\end{equation}
As explained in \cite{RW1}, 
the elliptic class of a Schubert variety $X_w$ for $w\in W$ is well-defined only after twisting by a fractional line bundle $\lambda\in \operatorname{Pic}(G/B)_{\mathbb{Q}}$ and the formal variable $\hbar$. 
In order to include the dependence on $\la$ and $\hbar$, 
we require the elliptic class of $X_w$ to have values in a localization of 
\begin{equation}\label{eq:defofM}
\mathcal{M}=K_T(G/B)[\operatorname{Pic}(G/B)^\v\oplus\mathbb{Z}\hbar][[q]],
\end{equation}
where $\operatorname{Pic}(G/B)^\v\cong H_2(G/B,\mathbb{Z})$ is the dual lattice of $\operatorname{Pic}(G/B)$.  
By identifying $(G/B)^T=W$, $K_T(\pt)=\mathbb{Q}[X^*(T)]$ and $H_2(G/B;\mathbb{Z})$ with the coroot lattice $X_*(T)$, we can embed $\mathcal{M}\subset \mathcal{Q}_W^*$ by localization.
As a generalized oriented cohomology theory, 
the identification \eqref{eq:defEGB} is a ring isomorphism, while the Poincar\'e pairings are different. 
We translate to the pairing over $\mathcal{Q}_W^*$ by 
\begin{equation}\label{eq:PoincarePair}
\langle f,g\rangle=
\langle f,g\rangle_{G/B}:=
\sum_{u\in W}
\frac{f(u)g(u)}{\prod_{\al>0}\theta(-z_{u\al})}\in \mathcal{Q}
\text{ where }f,g\in \mathcal{Q}_W^*.
\end{equation}

\subsection{Elliptic Schubert classes}

Recall the definition of $\mathsf{P}$ and $\mathsf{Q}$ in \eqref{eq:defofPQ} and the convention \eqref{eq:convtheta}.


\begin{Def}\label{def:EDL}
For a simple root $\alpha\in X^*(T)$, we define 
the \emph{elliptic Demazure-Lusztig (DL) operator}
\begin{equation}\label{eq:defEDL}
T_\al = \de_{\al}^{\dyn} \big(\P(z_\al,\la_\alv)+\Q(z_\al,\la_\alv)
\de_{\al}\big)\in \mathcal{Q}_{W^2},
\end{equation}
where 
$\de_{\al}=\de_{s_\al}$ and $\de^{\dyn}_{\al}=\de^{\dyn}_{s_\al}$. 
\end{Def}

\begin{Rmk}{\rm 
In \cite[Theorem 1.3]{RW1}, Rim\'anyi and Weber introduced the following operator for a simple root $\alpha$:
$$T^{RW}_\alpha = \left(\frac{\delta(z_\alpha,\lambda_{\alv})}{\delta(-\la_{\al^\v},h)}-
\frac{\delta(z_\alpha,h)}{\delta(-\la_{\al^\v},h)}\delta_{\al}\right)\delta^\dyn_{\al}
\;\;\;\;\;\text{ where }
\delta(a,b)=\frac{\theta(a+b)\theta'(0)}{\theta(a)\theta(b)}.$$
In other words, we have}
\begin{align*}
T_\al^{RW}&=\left(\four{z_\al+\la_\alv}{h}{z_\al}{\la_\alv-h}
-\four{z_\al+h}{\la_\alv}{z_\al}{\la_\alv-h}\delta_\al\right)\delta_\al^\dyn\\
& = \delta_\al^\dyn
\left(\four{z_\al-\la_\alv}{h}{z_\al}{-\la_\alv-h}
-\four{z_\al+h}{-\la_\alv}{z_\al}{-\la_\alv-h}\delta_\al\right)\\
& = \delta_\al^\dyn
\left(-\four{z_\al-\la_\alv}{h}{\la_\alv+h}{z_\al}
-\four{z_\al+h}{\la_\alv}{\la_\alv+h}{z_\al}\delta_\al\right) = -T_\alpha.
\end{align*}
{\rm Since $T_\alpha^{RW}$ satisfies the Weyl group relations, so does $T_\alpha$. 
In particular, $T_w$ is well-defined for any $w\in W$. }
\end{Rmk}

For $u\in W$, it is a standard fact about the twisted group algebra $\calQ_{W^2}$ that  one can expand 
\begin{align}\label{eq:defawubwu}
\de^{\dyn}_{u^{\-1}}T_u& 
= \sum_{w\leq u} a_{u,w}\cdot \de_w, &
\de_u & = \sum_{w\leq u} b_{u,w}\cdot 
\de^{\dyn}_{w^{\-1}}T_w, & 
\text{where }
a_{u,w},b_{u,w}\in \mathcal{Q}.
\end{align}
By definition, the matrix $(a_{u,w})$ is triangular with respect to the Bruhat order, and the matrix $(b_{u,w})$ is its inverse. 

\begin{Def}\label{def:Schubertcls}
For $w\in W$, we define the elliptic Schubert class
$\bfE_w\in \mathcal{Q}_{W}^*$ by 
$$\bfE_w(u)=b_{u,w}\in \mathcal{Q},$$
where $b_{w,u}$ is defined in \eqref{eq:defawubwu}. 
\end{Def}

This definition is the elliptic analogue of Schubert classes defined by Kostant-Kumar \cite{KK86}. 
On another hand, these classes form the dual basis of a certain normalization of the elliptic Schubert classes introduced in \cite{RW1}; see Theorem \ref{thm:dualbasis}. 

\begin{Eg}\label{eg:abuwlleq1}
{\rm We compute $a_{u,w}$ and $b_{u,w}$ for $\ell(u)\leq 1$. 
By definition, we have 
$b_{\id,w}=a_{\id,w}=\delta_{w,\id}^{\mathsf{Kr}}$. 
Let $\alpha$ be a simple root. Since 
$$
\de_{\al}^\dyn T_\al = \P(z_\al,\la_\alv)+\Q(z_\al,\la_\alv)\de_\al, $$
using \eqref{eq:PQinv}, we obtain
$$\de_{\al} = \P(\la_\alv,z_\al)+\Q(\la_\alv,z_\al)\de_{\al}^\dyn T_\al.$$
So $a_{s_\al,w}=b_{a_\al,w}=0$ unless $w\in \{\id,s_\al\}$, and 
\begin{align*}
a_{s_\al,\id} & = \P(z_\al,\la_\alv),&
a_{s_\al,s_\al} & = \Q(z_\al,\la_\alv),\\
b_{s_\al,\id} & = \P(\la_\alv,z_\al),&
b_{s_\al,s_\al} & = \Q(\la_\alv,z_\al).
\end{align*}
Note that one can switch between the $a$-coefficients and the $b$-coefficients by exchanging the variables. 
This is not a coincidence, as seen in Theorem {\rm \ref{thm:LZZ}} below. }
\end{Eg}

\begin{Eg}\label{eg:auubuu}
{\rm We compute $a_{u,w}$ and $b_{u,w}$ when $w=u$. 
It is easy to see that
$$T_u \in  \de_u^\dyn\left(b_{u,u}\de_u+\sum_{v<u} \mathcal{Q}\cdot \de_v\right).$$
If $s_\al u>u$, we have 
\begin{align*}
T_{s_\alpha u} 
& = T_\alpha T_u 
=\de_\al^\dyn(\P(z_\al,\la_\alv)+\Q(z_\al,\la_\alv)\de_\al)\cdot T_u\\
& \in 
\de^\dyn_{s_\alpha u}
\left(\Q(z_\al,\la_{u^{\-1}\alv})\cdot {}^{s_\al}b_{u,u}\cdot \de_{s_\alpha u}+\sum_{v<s_\al u}\mathcal{Q}\cdot \de_v\right).
\end{align*}
By induction, it is not hard to see that
$$a_{u,u}=
\prod_{\begin{subarray}{c}
\alpha>0\\
u\alpha<0
\end{subarray}}
\Q(z_{\-u\al},\la_\alv).$$
Since $b_{u,u}=a_{u,u}^{-1}$, by \eqref{eq:PQinv}, we have 
$$b_{u,u}=
\prod_{\begin{subarray}{c}
\alpha>0\\
u\alpha<0
\end{subarray}}
\Q(\la_\alv,z_{\-u\al}).$$}
\end{Eg}

\subsection{3D mirror symmetry}
We can repeat the construction for the Langlands dual group $G^L$ in the same twisted group algebra as follows. 
By analogy with Definition \ref{def:EDL}, we define 
\begin{equation}\label{eq:defdEDL}
T^{\dyn}_\al = \de_{\al} \big(\P(\la_\alv,z_\al)+\Q(\la_\alv,z_\al)
\de^{\dyn}_{\al}\big)\in \mathcal{Q}_{W^2}. 
\end{equation}
Clearly $T^{\dyn}_\alpha$ is obtained from $T_\alpha$ by switching the roles of $\la_\alv$ and $z_\al$, so it is the elliptic DL operator for the Langlands dual system. Similar as \eqref{eq:defawubwu}, we define 
\begin{align}\label{eq:defdawubwu}
\de_{u^{\-1}}T^{\dyn}_u& 
= \sum_{w\leq u} a^{\dyn}_{u,w}\cdot \de^{\dyn}_w,
&
\de^\dyn_u & = \sum_{w\leq u} b^\dyn_{u,w}\cdot 
\de_{w^{\-1}}T^{\dyn}_w, & 
\text{where }
a^{\dyn}_{u,w},b^{\dyn}_{u,w}\in \mathcal{Q}.
\end{align}
The matrix $(b^{\dyn}_{u,w})$ is the inverse of $(a^{\dyn}_{u,w})$. It turns out that 
$a^{\dyn}_{u,w},b^{\dyn}_{u,w}$ and $a_{u,w},b_{u,w}$ are related by the following simple relation that can be understood as another form of  \emph{3D mirror symmetry}.

\begin{Th}\cite{RW2, LZZ}\label{thm:LZZ} 
For any $u,w\in W$, we have 
$$a_{u^{\-1},w^{\-1}}=b_{u,w}^{\dyn},\qquad 
b_{u^{\-1},w^{\-1}}=a_{u,w}^{\dyn}.$$
\end{Th}
\begin{proof}
Since we are using a slightly different convention,  we include the proof for the reader's convenience. It is essentially the same as in \cite{LZZ}. 

If $\ell(u)\leq 1$, this is done in Example \ref{eg:abuwlleq1}. 
So it suffices to show that $a_{u,w}$ and $b_{u,w}^{\dyn}$ have the same recursion formula, and similarly for $a^{\dyn}_{u,w}$ and $b_{u,w}$. 
We first derive the recursion for $b_{u,w}^\dyn$. 
Starting from 
\[
\delta_{u} ^\dyn
=\sum_{w\in W} b_{u,w}^{\dyn}\cdot \de_{w^{\-1}} T_{w}^\dyn,
\]
the left multiplication by $\delta_\alpha^\dyn$ gives 
\begin{align*}
\delta_{s_\al u}^\dyn &=\de_\al^\dyn\sum_{w}b_{u,w}^{\dyn}\cdot \de_{w^{\-1}} T_w^\dyn
 =\sum_{w}{}^{s^\dyn_\al}b_{u,w}^{\dyn}\cdot \de_{w^{\-1}} \de_\al^\dyn T_w^\dyn\\
&=\sum_{w} {}^{s^\dyn_\al}b_{u,w}^{\dyn}\cdot \de_{w^{\-1}}(\P( z_\al, \la_{\alv})+\Q(z_\al, \la_{\alv})\de_\al T_\al^\dyn)T^\dyn_{w}\\
&=\sum_{w}
    {}^{s^\dyn_\al}b_{u,w}^{\dyn} \cdot 
    \P(z_{w^{\-1}\al}, \la_{\alv})\cdot 
    \de_{w^{\-1}}T_w^\dyn+
    {}^{s^\dyn_\al}b_{u,w}^{\dyn} \cdot 
    \Q( z_{w^{\-1}\al}, \la_\alv)\cdot
    \de_{w^{\-1}s_\al}T_{s_\al w}^\dyn\\
& =\sum_{w}\left(
    {}^{s^\dyn_\al}b_{u,w}^{\dyn}\cdot 
    \P( z_{w^{\-1}\al}, \la_{\alv})+
    {}^{s^\dyn_\al}b_{u,s_\al w}^{\dyn}\cdot
    \Q(\-z_{w^{\-1}\al}, \la_{\alv})\right)\de_{w^{-1}}T_w^\dyn.
\end{align*}
By comparing the coefficients, we  conclude that 
\begin{equation}\label{eq:rec-bdyn}
b_{s_\al u,w}^\dyn = 
    {}^{s^\dyn_\al}b_{u,w}^{\dyn}\cdot 
    \P( z_{w^{\-1}\al}, \la_{\alv})+
    {}^{s^\dyn_\al}b_{u,s_\al w}^{\dyn}\cdot
    \Q(\-z_{w^{\-1}\al}, \la_{\alv}).
\end{equation}

Now we consider  $a_{u,w}$. 
Starting from 
\[
\de_{u}^{\dyn} T_{u^{\-1}}
=\sum_{w\in W}  a_{u^{\-1},w^{\-1}}\cdot \de_{w^{\-1}},
\]
the right multiplication by $T_{\alpha}$ gives
\begin{align*}
T_{(s_{\al}u)^{\-1}} &=T_{u^{-1}}T_\al\\
& = \de^\dyn_{u^{\-1}}\sum_{w} a_{u^{\-1},w^{\-1}}\cdot \delta_{w^{\-1}}
\cdot \de^\dyn_\al
(\P(z_{\al},\la_{\al^\v})+\Q(z_{\al},\la_{\al^\v})\de_\al) \\
& = \de_{(s_\al u)^{\-1}}^\dyn\sum_{w} 
    {}^{s^\dyn_\al} {a}_{u^{\-1},w^{\-1}} \cdot
    \P(z_{w^{\-1}\al}, \la_{\alv})\cdot 
    \de_{w^{\-1}}+
    {}^{s^\dyn_\al} {a}_{u^{\-1},w^{\-1}} \cdot 
    \Q(z_{w^{\-1}\al}, \la_{\alv}))\cdot
    \de_{w^{\-1}s_\al}\\
& = \de_{(s_\al u)^{\-1}}^\dyn\sum_{w}\left(
    {}^{s^\dyn_\al} {a}_{u^{\-1},w^{\-1}} \cdot 
    \P(z_{w^{\-1}\al}, \la_{\alv})
    +
    {}^{s^\dyn_\al} {a}_{u^{\-1},w^{\-1} s_\al} \cdot 
    \Q(\-z_{w^{-1}\al}, \la_{\alv}) 
    \right)\de_{w^{\-1}}.
\end{align*}
By comparing the coefficients, we  conclude that
\begin{equation}\label{eq:rec-a}
{a}_{u^{\-1}s_\al,w^{\-1}}=
    {}^{s^\dyn_\al} {a}_{u^{\-1},w^{\-1}} \cdot 
    \P(z_{w^{\-1}\al}, \la_{\alv})
    +
    {}^{s^\dyn_\al} {a}_{u^{\-1},w^{\-1}s_\al} \cdot 
    \Q(\-z_{w^{\-1}\al}, \la_{\alv}).
\end{equation}
Since \eqref{eq:rec-a} and \eqref{eq:rec-bdyn} have the same form, we  conclude that 
$b^{\dyn}_{s_\al u,w}=a_{s_\al u,w}$. The proof for the second identity follows similarly.
\end{proof}

\subsection{Recursions for elliptic Schubert classes}
With the help of Theorem \ref{thm:LZZ}, we are able to establish two recursion formulas for $b_{w,u}$, which will be used to study the elliptic Schubert classes. 

\begin{Prop}\label{prop:rec-b}
For any $u,w\in W$ and a simple root $\alpha$, we have 
\begin{align}\label{eq:rec-bleft}
b_{s_\al u,w} & = 
    \P(\la_{w^{\-1}\alv},z_{\al})\cdot
    {}^{s_\al}b_{u,w} +
    \Q(\-\la_{w^{\-1}\alv},z_{\al})\cdot 
    {}^{s_\al}b_{u,s_\al w},\\
    \label{eq:rec-bright}
b_{us_\al,w} & = 
    \P(\la_{\alv}, z_{u\al})\cdot b_{u,w}+
    \Q(\la_{\alv}, z_{u\al})\cdot{}^{s^\dyn_\al} b_{u,ws_\al}.
\end{align}
\end{Prop}
\begin{proof}
The proof of \eqref{eq:rec-bleft} is similar to  the proof of the recursion formula for $b^\dyn_{s_\al u, w}$ in  \eqref{eq:rec-bdyn}, so we skip it. 
In order to prove \eqref{eq:rec-bright}, we need to use the fact $a_{u^{\-1},w^{\-1}}^{\dyn}=b_{u,w}$, which was shown in Theorem \ref{thm:LZZ} above. 
Multiplying by $T_\al^\dyn=\de_\al(\P(\la_{\alv}, z_\al)+\Q(\la_{\alv}, z_\al)\de_\al^\dyn)$ on the left of the identity
\[
T_{u^{\-1}}^\dyn=\de_{u^{\-1}}\sum_{w\in W} a^{\dyn}_{u^{\-1},w^{\-1}}\de_{w^{-1}}^\dyn, 
\]
we get
\begin{align*}
T^\dyn_{s_\al u^{\-1}}
& =\sum_w
    \de_\al(\P(\la_{\al^\v}, z_\al)+
    \de_\al \Q(\la_{\al^\v}, z_\al)
    \de_\al^\dyn)\de_{u^{\-1}}
    \cdot a^{\dyn}_{u^{\-1},w^{\-1}}\cdot
    \de_{w^{\-1}}^\dyn\\
& =\de_{s_\al u^{\-1}}\sum_w 
    \P(\la_{\alv}, z_{u\al})
    a_{u^{\-1},w^{\-1}}^{\dyn}\cdot\de_{w^{\-1}}^\dyn+
    \Q(\la_{\alv}, z_{u\al})\cdot 
    {}^{s^\dyn_\al}a_{u^{\-1},w^{\-1}}^{\dyn}\cdot \de_{s_\al w^{\-1}}^\dyn\\
&= \de_{s_\al u^{\-1}}\sum_w \left(
    \P(\la_{\alv}, z_{u\al})a_{u^{\-1},w^{\-1}}^{\dyn}
    +\Q(\la_{\alv}, z_{u\al}) \cdot 
    {}^{s^\dyn_\al }a_{u^{\-1},s_\al w^{\-1}}^{\dyn}
    \right)\de_{w^{\-1}}^\dyn.
\end{align*}
By comparing the coefficients, \eqref{eq:rec-bright} follows. 
\end{proof}

\begin{Rmk}{\rm The recursion \eqref{eq:rec-bright} is equivalent to the recursive generation of the elliptic classes via the operators $T_\alpha$, see~\cite[Section~3.2]{Zho} for the formulation in the same setup as the one here. }
\end{Rmk}



\subsection{Elliptic Schubert classes for $G/P$}
In order to generalize the above setting to a parabolic subgroup, we need to assume that $G$ is simply connected, i.e., $X_*(T)$ is the coroot lattice $\check{Q}$. 
We denote 
\begin{equation}\label{eq:paraboliccoroots}
\check{Q}^{P} = \bigoplus_{\alpha\in \Sigma\setminus \Sigma_P}\mathbb{Z}\alpha^\v\subset 
\bigoplus_{\alpha\in \Sigma}\mathbb{Z}\alpha^\v
=\check{Q}=X_*(T).
\end{equation}
We denote 
\begin{equation}\label{eq:defofQcal_para}
\mathcal{Q}^P=
\operatorname{Frac}\mathbb{Q}[\Lambda^P][[q]]
\text{ where }
\Lambda^P = X^*(T)\oplus \check{Q}^P\oplus
\mathbb{Z}\hbar,
\end{equation}
and 
\begin{equation}
\label{eq:defofQcalWstar_para}
\mathcal{Q}_{W/W_P}^*=\operatorname{Map}(W/W_P,\mathcal{Q}^P)
=\{f\in \operatorname{Map}(W,\mathcal{Q}^P):
\forall v\in W_P,\,f(uv)=f(u)\}.
\end{equation}

We now define the parabolic version of elliptic Schubert classes. 
We have the following linear map of lattices
\begin{equation}
\label{eq:pro}\check{Q}\oplus \mathbb{Z}\hbar
\twoheadrightarrow 
\check{Q}^P\oplus \mathbb{Z}\hbar,\qquad 
\alv\longmapsto 
\begin{cases}
\alv, & \al\notin \Sigma_P,\\
-\hbar, & \al\in \Sigma_P.
\end{cases}
\end{equation}
Denote by $\mathcal{Q}'\subset \mathcal{Q}$  the subset of elements having no poles along $\la_\alv=-\hbar$ for all  roots $\al\in \Phi$. 
Then the map \eqref{eq:pro} induces a map 
\begin{equation}\label{eq:spelato-h}
   \mathcal{Q}'\longrightarrow \mathcal{Q}^P,\qquad c\longmapsto [c]_P,
\end{equation}
by specializing $\la_{\alv}$ to $-\hbar$ for all $\al\in \Sigma_P$. 

We need the following classical result; see, e.g., \cite[Cor.~2.5.2]{BB}.

\begin{Lemma}\label{lem:Deodhar}
Let  $w\in W^P$ and  $\sigma$ be a simple reflection. Then either $\sigma w\in W^P$ or $\sigma w =w \tau>w$  for a simple reflection $\tau\in W_P$. 
\end{Lemma}

\begin{Lemma}\label{lem:buw_para}
For $u,w\in W$, we have 
\begin{enumerate}[\quad \rm(1)]
    \item \label{lem:buw_para1}
    $b_{u,w}\in \mathcal{Q}'$;
    \item \label{lem:buw_para2}
    $[b_{u,w}]_P=[b_{uv,w}]_P$ for any $v\in W_P$; 
    \item \label{lem:buw_para3}
    $[b_{u,w}]_P = 0$ for any  $w\not\in W^P$. 
\end{enumerate}
\end{Lemma}
\begin{proof}
\eqref{lem:buw_para1}
It is easy to see that $\calQ'$ is $W\times W^\dyn$-invariant. 
Note that $\P(\la_\alv,z_\al)$ and $\Q(\la_\alv,z_\al)$ only have poles along $\{\la_\alv=0\}$ and $\{z_\al=-\hbar\}$, so the conclusion follows directly by induction and the recursion \eqref{eq:rec-bleft}.

\eqref{lem:buw_para2}
It suffices to prove the identity for $v=s_\al$ with  $\al\in \Sigma_P$.
By \eqref{lem:buw_para1}, we can apply $[\_]_P$ on each term of \eqref{eq:rec-bright}. Since $[\P(\la_{\al^\v}, z_{u\al})]_P=1$ and $[\Q(\la_{\al^\v}, z_{u\al})]_P=0$, we obtain $[b_{us_\al,w}]_P = [b_{u,w}]_P$. 

\eqref{lem:buw_para3}
We prove this fact by induction on $u$. 
It is obvious when $u=\id$. 
Assuming it holds for $u$, we show it also holds for $s_\al u>u$ for any simple root $\al\in \Sigma$. 
By \eqref{lem:buw_para1}, we can apply $[\_]_P$ on each term of \eqref{eq:rec-bleft}. Let $w\not\in W^P$. There are two cases:
\begin{enumerate}[(a)]
\item 
If $s_\alpha w\notin W^P$, since $[{}^{s_\al}b_{u,w}]_P={}^{s_\al}[b_{u,w}]_P=0$,  then \eqref{eq:rec-bleft} implies directly that $[b_{s_\al u,w}]_P=0$.
\item 
If $s_{\alpha}w\in W^P$, we can apply  Lemma \ref{lem:Deodhar} to $s_\alpha w$, and obtain that $\be:=-w^{-1}\al\in \Sigma_P^+$. So  $\be^\v=-w^{\-1}\alv$ is a simple coroot of $P$, which implies that $[\Q(-\la_{w^{-1}\alv}, z_\al)]_P=0$. So \eqref{eq:rec-bleft} again implies \mbox{$[b_{s_\al u,w}]_P=[{}^{s_\al}b_{u,w}]_P=0$}.
\qedhere
\end{enumerate}
\end{proof}

\begin{Def}\label{def:paraSchubertcls}
For $w\in W^P$, we define the \emph{parabolic elliptic Schubert class}
$\bfE_w^P\in \mathcal{Q}_{W/W_P}^*$ by
$$\bfE_w^P(u)=[b_{u,w}]_P\in \mathcal{Q}^P.$$
By Lemma {\rm \ref{lem:buw_para}}, the class $\bfE_w^P$ is well defined. 
\end{Def}

See Corollary \ref{coro:paradualbasis} for the relation to the parabolic elliptic classes defined in \cite{KRW}.

\begin{Rmk}{\rm 
We note that the properties in Lemma \ref{lem:buw_para} only hold for elliptic cohomology. This phenomenon ceases upon degeneration to $K$-theory (i.e. to Segre motivic classes) since the degeneration and the specialization are not compatible. }
\end{Rmk}

\section{Elliptic Billey Formula}
\label{sec:ellBilfor}

In this section, we give a Billey-type formula for the elliptic classes $\bfE_w$.  
In type $A$, it can be represented by a diagrammatic formula. 

\subsection{The dynamical $R$-matrix}
We consider
$$\mathcal{Q}_{W^\dyn} = \mathcal{Q}\rtimes \mathbb{Q}[W^\dyn].$$
It is a free left $\mathcal{Q}$-module with basis $\de_{v}^\dyn$, and there is a  multiplication given by 
$$
a\,\de_v^{\dyn} \cdot 
a'\,\de_{v'}^{\dyn}
= a\cdot {}^{v^\dyn}a'\,\de^\dyn_{vv'}. $$
For a simple root $\al_i$ and a weight $\beta\in X^*(T)$, we consider the \emph{dynamical $R$-matrix}
\begin{equation}\label{eq:dynRmat}
h_i(\beta) = \P(\la_{\al^\v_i},z_{\beta})+\de_{i}^\dyn\cdot \Q(\la_{\al^\v_i},z_{\beta})\in \mathcal{Q}_{W^\dyn},
\end{equation}
where $\de^\dyn_i = \de^\dyn_{s_i}$.
We are going to construct a generating function of $b_{u,w}=\bfE_w(u)$ via dynamical $R$-matrices.

\begin{Rmk} {\rm The dynamical $R$-matrices have appeared in Felder's work on elliptic quantum groups \cite{F95}.}
\end{Rmk}

We first fix a (not necessarily reduced) word of $u\in W$
\begin{equation}\label{eq:decomofu}
u=s_{i_1} \cdots  s_{i_\ell}, 
\end{equation}
which determines a sequence of roots:
$$\beta_j = s_{i_1}\cdots s_{i_{j-1}}\al_{i_j}, ~j=1,\ldots,\ell. 
$$
If  the decomposition \eqref{eq:decomofu} is reduced, $\beta_j$ are all positive, and 
we have 
\begin{align*}
\{\beta_1,\ldots,\beta_\ell\}
& =\{\beta>0: w^{-1}\beta<0\}.
\end{align*}
Our first main result is the following theorem.

\begin{Th}\label{th:BilleyAsGF}
With the notation introduced above, we have 
\begin{equation}\label{eq:BilleyAsGF}
h_{i_1}(\beta_1)\cdots h_{i_\ell}(\beta_\ell)
=\sum_{w\in W} \de^\dyn_{w}\cdot \bfE_w(u).
\end{equation}
In particular, the left-hand side does not depend on the choice of decomposition \eqref{eq:decomofu}. 
\end{Th}
\begin{proof}
We use induction on the length $\ell$ of $u$. 
The base case $u=\id$ is obvious. 
Assuming that the theorem is true for all $u$ with length $\ell$, we prove it for $us_{\al}$ with length $\ell+1$ and $\al$  a simple root. In this case $\beta_{\ell+1}=u\al$. Recall that $b_{u,w}=\bfE_w(u)$. By induction, we have
\begin{align*}
h_{i_1}(\beta_1)\cdots h_{i_\ell}(\beta_\ell)\cdot h_{i_{\ell+1}}(\beta_{\ell+1})
& = 
\left(
\sum_{w\in W}\de^\dyn_w\cdot b_{u,w}\right)
\left(\P(\la_\alv,z_{u\al})+\de_{\al}^\dyn\cdot \Q(\la_\alv,z_{u\al})\right)\\
& =
\sum_{w\in W}
\de^\dyn_w\cdot b_{u,w}\cdot \P(\la_\alv,z_{u\al})
+
\de^\dyn_{ws_\al}
\cdot {}^{s_\al^\dyn}b_{u,w}\cdot \Q(\la_\alv,z_{u\al})\\
& = 
\sum_{w\in W}
\de^\dyn_w\left(
b_{u,w}\cdot \P(\la_\alv,z_{u\al})+
{}^{s_\al^\dyn}b_{u,w s_\al }\cdot \Q(\la_\alv,z_{u\al})
\right)\\
& = \sum_{w\in W}
\de^\dyn_w\cdot b_{us_\al,w}\qquad 
\text{ by \eqref{eq:rec-bright}.}\qedhere
\end{align*}
\end{proof}

\begin{Coro}\label{coro:YBequation}
The dynamical $R$-matrices \eqref{eq:dynRmat} satisfy the following properties. 
\begin{itemize}
    \item Unitary property
    $$h_{i}(x)h_{i}(-x)=1.$$
    \item Yang-Baxter equations
    \begin{align*}
        h_i(x)h_j(y)&=h_j(y)h_i(x),&\text{ if }(s_is_j)^2=\id;\\
        h_i(x)h_j(x+y)h_i(y)&=h_j(y)h_i(x+y)h_j(x),&
        \text{ if }(s_is_j)^3=\id .
    \end{align*}
\end{itemize}
There are similar expressions for $(s_is_j)^4=\id$ and  $(s_is_j)^6=\id$, see also \cite[Prop 3.1]{billey}.
\end{Coro}
\begin{proof}
We first consider the root system $\mathsf{A}_2$.
$$
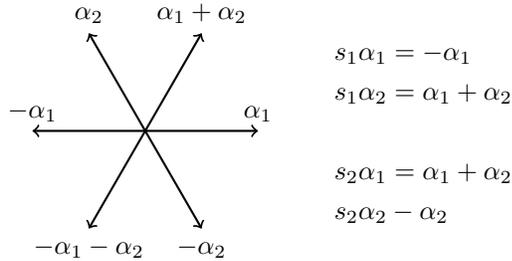
\begin{matrix}\begin{tikzpicture}[scale = 0.75]
\draw[->,thick] (0,0) --  (2,0) node[above]{$\al_1$};
\draw[->,thick] (0,0) to (-2,0) node[above]{$-\al_1$};
\draw[->,thick] (0,0) to (1,1.732) node[above]{$\al_1+\al_2$};
\draw[->,thick] (0,0) to (1,-1.732)node[below]{$-\al_2$};
\draw[->,thick] (0,0) to (-1,1.732)node[above]{$\al_2$};
\draw[->,thick] (0,0) to (-1,-1.732)node[below]{$-\al_1-\al_2$};
\end{tikzpicture}
\end{matrix}\qquad 
\begin{aligned}
s_1\alpha_1 &=-\alpha_1 \\
s_1\alpha_2 &= \alpha_1+\alpha_2 \\\\
    s_2\alpha_1 & = \alpha_1+\alpha_2\\ 
    s_2\alpha_2 & - \alpha_2
\end{aligned}
$$
Consider the two decompositions 
$$
u=s_1s_2s_1\quad \text{ vs. }\quad
u=s_2s_1s_2.$$
By Theorem \ref{th:BilleyAsGF}, we have 
$$
h_1(\al_1)h_2(\al_1+\al_2)h_1(\al_2)
=h_2(\al_2)h_1(\al_1+\al_2)h_2(\al_1).$$
By replacing the two independent weights $\al_1,\al_2$ with the symbols $x$ and $y$, we get the Yang-Baxter equation for $(s_is_j)^3=\id$.

Similarly, the other Yang-Baxter equations can be obtained by considering the other root systems of rank $2$. The unitary property can be proved via the root system of rank $1$ using $\id=s^2$.
\end{proof}

\subsection{Billey-type formula}
Now we  expand the formula in Theorem \ref{th:BilleyAsGF}. Recall that we fixed a decomposition \eqref{eq:decomofu}. 
For a subset $J\subseteq \{1,\ldots,\ell\}$, we introduce 
$\epsilon_{j}=\delta_{j\in J}^{\mathsf{Kr}}\in \{0,1\}$ for $1\leq j\leq \ell$, and define the subword of \eqref{eq:decomofu}
$$w(J)= s_{i_1}^{\epsilon_1}\cdots s_{i_\ell}^{\epsilon_\ell}=
\prod_{j\in J}^{\rightarrow}s_{i_j}\leq u.$$
For $j=1,\ldots,\ell$, we consider the coroots:
\begin{equation}\label{eq:defofgamma}
\check{\gamma}_{j}^{J} = s_{i_{\ell}}^{\epsilon_{\ell}}
\cdots 
s_{i_{j+1}}^{\epsilon_{j+1}}
\alpha^\v_{i_j}.
\end{equation}

\begin{Eg}{\rm 
Consider the  root system $A_3$. 
We illustrate the formula for 
$$u=s_1s_2s_3s_2s_1s_3,\qquad J = \{2,4,5\}.$$
Then 
$$w = \cancel{s_1}s_2\cancel{s_3}s_2s_1\cancel{s_3}=s_2s_2s_1=s_1.$$
By definition, 
$$\begin{aligned}
\beta_1 & = \alpha_1,&&\qquad&
    \check{\gamma}^J_1 & = s_1s_2s_2\alpha^\v_1,&&\qquad&
    1\notin J,\\
\beta_2 & = s_1\alpha_2,&&\qquad&
    \check{\gamma}^J_2 & = s_1s_2\alpha^\v_2 ,&&\qquad&
    2 \in J,\\
\beta_3 & = s_1s_2\alpha_3,&&\qquad&
    \check{\gamma}^J_3 & = s_1s_2\alpha^\v_3,&&\qquad&
    3\notin J,\\
\beta_4 & = s_1s_2s_3\alpha_2,&&\qquad&
    \check{\gamma}^J_4 & = s_1\alpha^\v_2,&&\qquad&
    4\in J,\\
\beta_5 & = s_1s_2s_3s_2\alpha_1,&&\qquad&
    \check{\gamma}^J_5 & = \alpha^\v_1,&&\qquad&
    5\in J,\\
\beta_6 & = s_1s_2s_3s_2s_1\alpha_3,&&\qquad&
    \check{\gamma}^J_6 & = \alpha^\v_3,&&\qquad&
    6\notin J.
\end{aligned}
$$}
\end{Eg}

\begin{Th}\label{th:Billey}
With the notation introduced above, 
we have
$$\bfE_w(u)=
\sum_{w(J)=w}\prod_{j=1}^\ell\begin{cases}
\Q(\la_{\check{\gamma}_{j}^{J}},z_{\beta_j}),& j\in J,\\
\P(\la_{\check{\gamma}_{j}^{J}},z_{\beta_j}),& j\notin J.
\end{cases}$$
\end{Th}
\begin{proof}
We can expand \eqref{eq:BilleyAsGF} as 
$$\sum_{J}\prod_{j=1}^{\ell}
\begin{cases}
\de^\dyn_{i_j}\Q(\la_{\al^\v_{i_j}},z_{\beta_j}), & j\in J,\\
\P(\la_{\al^\v_{i_j}},z_{\beta_j}), & j\notin J.
\end{cases}
=\sum_{J}\de^\dyn_{w(J)}
\prod_{j=1}^\ell
\begin{cases}
\Q(\la_{\check{\gamma}_{j}^{J}},z_{\beta_j}), & j\in J,\\
\P(\la_{\check{\gamma}_{j}^{J}},z_{\beta_j}), & j\notin J.
\end{cases}$$
By comparing the coefficients, the formula follows.
\end{proof}
We would like to emphasize that the set of coroots $\check{\gamma_j}^J$ depends not only on $u$, but also on $J$. 
Now we  state the parabolic version of Theorem \ref{th:Billey}. 
In fact, it gives another proof of \eqref{lem:buw_para3} in Lemma~\ref{lem:buw_para}.

\begin{Th}\label{th:Billeypara}
Using the notation introduced above, we have 
$$\bfE_w^P(u)=\sum_{w(J)=w}
\prod_{j=1}^\ell\begin{cases}
{}
[\Q(\la_{\check{\gamma}_{j}^{J}},z_{\beta_j})]_P,& j\in J,\\
{}[\P(\la_{\check{\gamma}_{j}^{J}},z_{\beta_j})]_P,& j\notin J,
\end{cases}$$
where the sum over $J\subset \{1,\ldots,\ell\}$ such that 
\begin{equation}\label{eq:paracondition}
\forall j =1,\ldots,\ell, \text{we have }s_{i_j}^{\epsilon_j}
\cdots 
s_{i_\ell}^{\epsilon_\ell}\in W^P. 
\end{equation}
In particular, $\bfE_w^P(u)=0$ unless $w\in W^P$. 
\end{Th}
\begin{proof}
By Definition \ref{def:paraSchubertcls}, the element $\bfE^P_w(u)$ can be obtained directly from $\bfE_w(u)$ via specialization.  Note that 
\[
[\P(\la_{\check{\gamma}^J_{j}} , z_{\be_j})]_P=1,~ [\Q(\la_{\check{\gamma}^J_{j}} , z_{\be_j})]_P=0, ~\text{ if }\check{\gamma}^J_j\in \Sigma_P^\v. 
\]
Fix a $J=\{j_1<j_2<\cdots <j_k\}$, and denote \[
w_{j}=s_{i_{j}}^{\epsilon_j}\cdots s_{i_{\ell}}^{\epsilon_{\ell}}. 
\]
Note $w_{j_k}=s_{i_{j_k}}$ is a simple reflection with $\check{\gamma}_{{j_k}}^J=\al_{i_{j_k}}^\v$. 
To prove the theorem, it suffices to show that for $J$ such that $w_{j_t}\not\in W^P$, $\check{\gamma}^J_{j_t}\in \Sigma^\v_P$. 

If $w_{j_k}\not\in W^P$, then $\check{\gamma}_{j_k}^J=\al^\v_{i_{j_k}}\in \Sigma_P^\v$. 
Now let $w_{j_k}\in W^P$, and let $t$ be the maximal number so that $w_{j_t}\not\in W^P$ and  $w_{j_{t+1}}\in W^P$. Then $w_{j_t}=s_{i_{j_t}}w_{j_{t+1}}\not\in W^P$ implies that $w_{j_t}=w_{j_{t+1}}s_\be$ for some $\be\in \Sigma_P$. Thus, 
\[
\check{\gamma}_{j_t}^J=w_{j_{t+1}}^{-1}\al_{i_{j_t}}^\v=\be^\vee\in \Sigma_P^\v.
\]
The proof is finished. 
\end{proof}

\section{Diagrams for type A}
\label{sec:diatypeA}

Now we switch to $\operatorname{GL}_n$, 
where we can identify 
$$X^*(T) = \mathbb{Z}\epsilon_1\oplus \cdots \oplus \mathbb{Z}{\epsilon_n}=X_*(T)$$
with simple roots $\epsilon_i-\epsilon_{i+1}$, for $i=1,\ldots,n-1$. 
We denote 
$$z_{i}=z_{\epsilon_i}\in \mathcal{Q},\qquad \la_i=\la_{\epsilon_i}\in \mathcal{Q}.$$
The Weyl group $W=S_n$ is the symmetric group. 

\subsection{Wiring diagrams}
Let $u$ be a permutation. A choice of a decomposition \eqref{eq:decomofu} can be illustrated by a \emph{wiring diagram}, where each $s_{i_j}$ is represented by a crossing of  strings of heights $i_j$ and $i_j+1$ in  the $j$-th interval. 
A choice of $J\subseteq \{1,\ldots,\ell\}$, corresponding to a choice of a subword, can be represented by a ``sub-wiring diagram'' obtained by resolving $s_{i_j}$ for $j\notin J$. 
To distinguish them, we color the strings in the wiring diagram of $u$ by {\color{red}red}, and those in the sub-wiring diagram by {\color{blue}blue}. 
We label the {\color{red}red} (resp., {\color{blue}blue}) strings by their {\color{red}targets} (resp., {\color{blue}sources}). 
Here is an example. 
$$
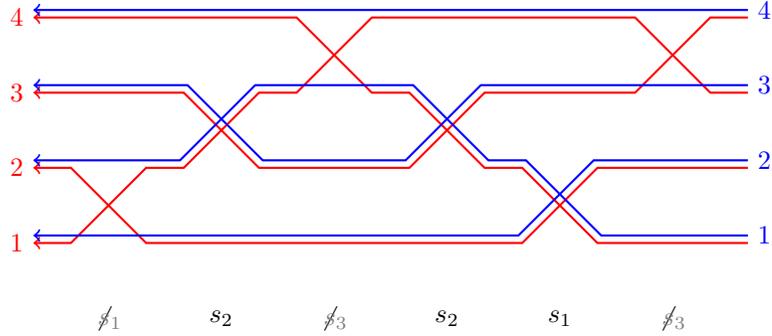
\begin{tikzpicture}
\redstrings
\draw[red,thick,<-] (0,1)node [left] {1}
    \b\u\u\u\h\h\d\e; 
\draw[red,thick,<-] (0,2)node [left] {2}
    \b\d\h\h\h\u\h\e; 
\draw[red,thick,<-] (0,3)node [left] {3}
    \b\h\d\h\u\h\u\e; 
\draw[red,thick,<-] (0,4)node [left] {4}
    \b\h\h\d\d\d\h\e; 
\draw[draw=none]
(0,0)--++(.5,0)--node{$\cancel{s_1}$}
    ++(1,0)--++(.5,0)--node{$s_2$}
    ++(1,0)--++(.5,0)--node{$\cancel{s_3}$}
    ++(1,0)--++(.5,0)--node{$s_2$}
    ++(1,0)--++(.5,0)--node{$s_1$}
    ++(1,0)--++(.5,0)--node{$\cancel{s_3}$}
    ++(1,0)--++(.5,0);
\bluestrings
\draw[blue,thick,<-](0,1.1)
    \b\h\h\h\h\u\h\e node [right] {2}; 
\draw[blue,thick,<-] (0,2.1)
    \b\h\u\h\d\d\h\e node [right] {1};
\draw[blue,thick,<-] (0,3.1) 
    \b\h\d\h\u\h\h\e node [right] {3};
\draw[blue,thick,<-] (0,4.1)
    \b\h\h\h\h\h\h\e node [right] {4};
\end{tikzpicture}$$
Then we can read $\beta_{j}$ and $\check{\gamma}^J_j$ from the sub-wiring diagram. At the $j$-th crossing, we have the following possibilities.
$$
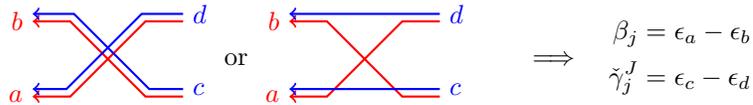
\begin{matrix}
\begin{tikzpicture}
\redstrings
\draw[red,thick,<-] (0,0) node [left] {$a$}
    \b\u\e;
\draw[red,thick,<-] (0,1) node [left] {$b$}
    \b\d\e;
\bluestrings
\draw[blue,thick,<-] (0,0.1) 
    \b\u\e node [right] {$d$};
\draw[blue,thick,<-] (0,1.1) 
    \b\d\e node [right] {$c$};
\end{tikzpicture}
\end{matrix}
\text{ or }
\begin{matrix}
\begin{tikzpicture}
\redstrings
\draw[red,thick,<-] (0,0) node [left] {$a$}
    \b\u\e;
\draw[red,thick,<-] (0,1) node [left] {$b$}
    \b\d\e;
\bluestrings
\draw[blue,thick,<-] (0,0.1) 
    \b\h\e node [right] {$c$};
\draw[blue,thick,<-] (0,1.1) 
    \b\h\e node [right] {$d$};
\end{tikzpicture}
\end{matrix}
\qquad \Longrightarrow\quad
\begin{aligned}
\beta_j  & = \epsilon_a-\epsilon_b\\
\check{\gamma}_j^J & = \epsilon_c-\epsilon_d
\end{aligned}$$
For instance, in the above example, we have:
$$\begin{aligned}
\beta_1 & = \epsilon_1-\epsilon_2,&&\qquad&
    \check{\gamma}_1^J & = \epsilon_2-\epsilon_1, \\
\beta_2 & = \epsilon_1-\epsilon_3,&&\qquad&
    \check{\gamma}_2^J & = \epsilon_3-\epsilon_1, \\
\beta_3 & = \epsilon_1-\epsilon_4,&&\qquad&
    \check{\gamma}_3^J & = \epsilon_1-\epsilon_4, \\
\beta_4 & = \epsilon_3-\epsilon_4,&&\qquad&
    \check{\gamma}_4^J & = \epsilon_1-\epsilon_3, \\
\beta_5 & = \epsilon_1-\epsilon_4,&&\qquad&
    \check{\gamma}_5^J & = \epsilon_1-\epsilon_2, \\
\beta_6 & = \epsilon_3-\epsilon_2,&&\qquad&
    \check{\gamma}_6^J & = \epsilon_3-\epsilon_4. 
\end{aligned}$$
This leads us to define the weight of a sub-wiring diagram to be the product of local weights below. 
\begin{equation}\label{eq:Ellweights}
\begin{matrix}
&\begin{tikzpicture}
\redstrings
\draw[red,thick,<-] (0,0) node [left] {$a$}
    \b\u\e;
\draw[red,thick,<-] (0,1) node [left] {$b$}
    \b\d\e;
\bluestrings
\draw[blue,thick,<-] (0,0.1) 
    \b\u\e node [right] {$d$};
\draw[blue,thick,<-] (0,1.1) 
    \b\d\e node [right] {$c$};
\end{tikzpicture}& 
\begin{tikzpicture}
\redstrings
\draw[red,thick,<-] (0,0) node [left] {$a$}
    \b\u\e;
\draw[red,thick,<-] (0,1) node [left] {$b$}
    \b\d\e;
\bluestrings
\draw[blue,thick,<-] (0,0.1) 
    \b\h\e node [right] {$c$};
\draw[blue,thick,<-] (0,1.1) 
    \b\h\e node [right] {$d$};
\end{tikzpicture}\\
\operatorname{weight}&
\Q(\la_{c}\-\la_d,z_a\-z_b)&
\P(\la_{c}\-\la_d,z_a\-z_b)
\end{matrix}
\end{equation}

Theorem \ref{th:Billey} can now be restated as follows. 

\begin{Th}\label{th:BilleytypeA}
Let $u,w\in S_n$. 
For a fixed wiring diagram $D_0$ of $u$, we have
$$\bfE_w(u)=\sum_{D}\operatorname{weight}(D),$$
where the sum is over all sub-wiring diagrams $D$ of $D_0$ whose permutation is $w$.
\end{Th}

\subsection{3D mirror symmetry again}
Next, we give a combinatorial explanation of Theorem \ref{thm:LZZ} in type $A$. 
Note that the theorem is equivalent to 
\begin{equation}\label{eq:3Dequivform}
\sum_{w\in W} b_{u,w}b^{\dyn}_{w^{\-1},v^{\-1}} = \delta^{\mathsf{Kr}}_{u,v},\qquad \forall u,v\in W.
\end{equation}
We prove \eqref{eq:3Dequivform} combinatorially.

First we remark that $b^{\dyn}_{w^{\-1},v^{\-1}}$ has a similar combinatorial formula to the one in Theorem \ref{th:BilleytypeA}, with the roles of $z_i$ and $\la_i$ switched.  
This leads to a description of 
the sum on the left-hand side of \eqref{eq:3Dequivform} via ``sub-sub-wiring diagrams,'' which are obtained by superimposing a sub-wiring diagram of a sub-wiring diagram over the original one. Here is an example.
$$
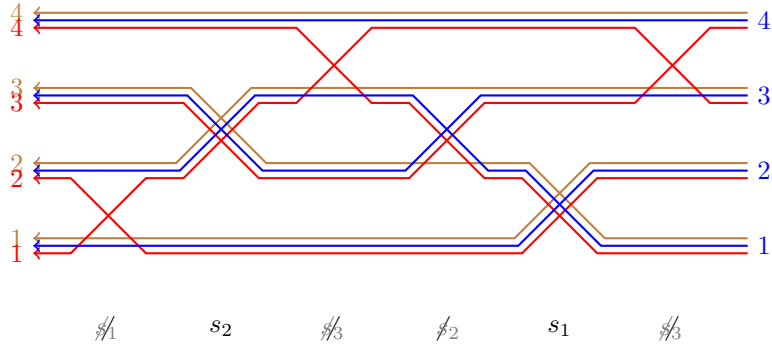
\begin{tikzpicture}
\def\b{--++(.25,0)}
\def\u{--++(.15,0)--++(1,1)--++(.35,0)}
\def\d{--++(.35,0)--++(1,-1)--++(.15,0)}
\def\h{--++(1.5,0)}
\def\e{--++(.25,0)}
\draw[brown,thick,<-] (0,1.2) node [left] {1} 
    \b\h\h\h\h\u\h\e;
\draw[brown,thick,<-] (0,2.2) node [left] {2}
    \b\h\u\h\h\h\h\e;
\draw[brown,thick,<-] (0,3.2) node [left] {3}
    \b\h\d\h\h\d\h\e;
\draw[brown,thick,<-] (0,4.2) node [left] {4}
    \b\h\h\h\h\h\h\e;
\redstrings
\draw[red,thick,<-] (0,1)node [left] {1}
    \b\u\u\u\h\h\d\e; 
\draw[red,thick,<-] (0,2)node [left] {2}
    \b\d\h\h\h\u\h\e; 
\draw[red,thick,<-] (0,3)node [left] {3}
    \b\h\d\h\u\h\u\e; 
\draw[red,thick,<-] (0,4)node [left] {4}
    \b\h\h\d\d\d\h\e; 
\draw[draw=none]
(0,0)--++(.5,0)--node{$\cancel{\!\!\cancel{s_1}}$}
    ++(1,0)--++(.5,0)--node{$s_2$}
    ++(1,0)--++(.5,0)--node{$\cancel{\!\!\cancel{s_3}}$}
    ++(1,0)--++(.5,0)--node{$\cancel{s_2}$}
    ++(1,0)--++(.5,0)--node{$s_1$}
    ++(1,0)--++(.5,0)--node{$\cancel{\!\!\cancel{s_3}}$}
    ++(1,0)--++(.5,0);
\bluestrings
\draw[blue,thick,<-](0,1.1)
    \b\h\h\h\h\u\h\e node [right] {2}; 
\draw[blue,thick,<-] (0,2.1)
    \b\h\u\h\d\d\h\e node [right] {1};
\draw[blue,thick,<-] (0,3.1) 
    \b\h\d\h\u\h\h\e node [right] {3};
\draw[blue,thick,<-] (0,4.1)
    \b\h\h\h\h\h\h\e node [right] {4};
\end{tikzpicture}$$
The weights are given as follows.
$$
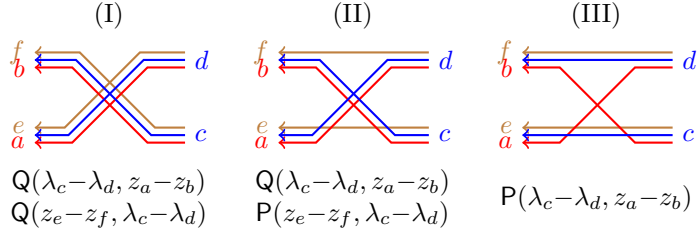
\begin{matrix}
\text{(I)}&\text{(II)}&\text{(III)}\\
\begin{tikzpicture}
\def\b{--++(.25,0)}
\def\u{--++(.15,0)--++(1,1)--++(.35,0)}
\def\d{--++(.35,0)--++(1,-1)--++(.15,0)}
\def\h{--++(1.5,0)}
\def\e{--++(.25,0)}
\draw[brown,thick,<-] (0,0.2) node [left] {$e$}
    \b\u\e;
\draw[brown,thick,<-] (0,1.2) node [left] {$f$}
    \b\d\e;
\redstrings
\draw[red,thick,<-] (0,0) node [left] {$a$}
    \b\u\e;
\draw[red,thick,<-] (0,1) node [left] {$b$}
    \b\d\e;
\bluestrings
\draw[blue,thick,<-] (0,0.1) 
    \b\u\e node [right] {$d$};
\draw[blue,thick,<-] (0,1.1) 
    \b\d\e node [right] {$c$};
\end{tikzpicture}& 
\begin{tikzpicture}
\draw[brown,thick,<-]
(0,0.2) node [left] {$e$}--++(.5,0)--
    ++(1,0)--++(.5,0);
\draw[brown,thick,<-]
(0,1.2) node [left] {$f$}--++(.5,0)--
    ++(1,0)--++(.5,0);
\redstrings
\draw[red,thick,<-] (0,0) node [left] {$a$}
    \b\u\e;
\draw[red,thick,<-] (0,1) node [left] {$b$}
    \b\d\e;
\bluestrings
\draw[blue,thick,<-] (0,0.1) 
    \b\u\e node [right] {$d$};
\draw[blue,thick,<-] (0,1.1) 
    \b\d\e node [right] {$c$};
\end{tikzpicture}
&
\begin{tikzpicture}
\draw[brown,thick,<-]
(0,0.2) node [left] {$e$}--++(.5,0)--
    ++(1,0)--++(.5,0);
\draw[brown,thick,<-]
(0,1.2) node [left] {$f$}--++(.5,0)--
    ++(1,0)--++(.5,0);
\redstrings
\draw[red,thick,<-] (0,0) node [left] {$a$}
    \b\u\e;
\draw[red,thick,<-] (0,1) node [left] {$b$}
    \b\d\e;
\bluestrings
\draw[blue,thick,<-] (0,0.1) 
    \b\h\e node [right] {$c$};
\draw[blue,thick,<-] (0,1.1) 
    \b\h\e node [right] {$d$};
\end{tikzpicture}\\
\begin{matrix}
\Q(\la_{c}\-\la_d,z_a\-z_b)\\
\Q(z_e\-z_f,\la_{c}\-\la_d)
\end{matrix}& 
\begin{matrix}
\Q(\la_{c}\-\la_d,z_a\-z_b)\\
\P(z_e\-z_f,\la_{c}\-\la_d)
\end{matrix}&
\P(\la_{c}\-\la_d,z_a\-z_b)
\end{matrix}$$
We can now explain the proof of \eqref{eq:3Dequivform}. 
We pick any reduced decomposition \eqref{eq:decomofu} of $u$. 
When $u=v$, there is only one sub-sub-wiring diagram, which is made up by (I) only. 
In particular, we can only have $e=a$ and $f=b$ in (I) above. By \eqref{eq:PQinv}, the weight is simply $1$. 
When $u\neq v$, we define an involution on sub-sub-wiring diagrams by changing the leftmost (II) (resp. (III)) to (III) (resp. (II)). 
Due to this choice (of the leftmost configuration), we still have $e=a$ and $f=b$. 
By \eqref{eq:PQinv}, the sum of their weights is zero. 
This gives a combinatorial proof of \eqref{eq:3Dequivform}.

\subsection{Parabolic version}
We now rephrase Theorem \ref{th:Billeypara} in the type $A$ case. 
%
%
Theoretically, we need to switch to $G=\operatorname{SL}_n$, but everything below can be lifted to $\operatorname{GL}_n$.

Let $\mathbf{a}=(a_1,\ldots,a_m)\in \mathbb{N}^m$ be a composition with 
$|\mathbf{a}|:=a_1+\cdots+a_m=n$. We have a parabolic subgroup $P\subset \operatorname{GL}_n$ whose Weyl group is the Young subgroup $W_P=S_{\mathbf{a}}:=S_{a_1}\times \cdots \times S_{a_m}\subset S_n$. 
We define 
$$A =\{a_1+\cdots+a_i:0\leq i\leq n\}= 
\{1,a_1+1,a_1+a_2+1,\ldots,a_1+\cdots+a_{m-1}+1\}.$$
Then the specialization map \eqref{eq:spelato-h} $f\mapsto [f]_P$ is given by 
$\la_{i}-\la_{i+1}\mapsto -\hbar$ for $i\notin A$. 
We denote $\la_a=[\la_a]_P$ for $a\in A$ by an abuse of notation. Then for any $i=1,\ldots,n$, 
$$\la_i\mapsto \la_{a}+(i-a)\hbar
\quad\text{where}\quad
a=\min(a\in A:a\leq i).$$
To compute $\bfE_w^P(u)$ using the sub-wiring diagrams introduced above, 
we relabel the {\color{blue}blue} strings by replacing the label $i$ by $a+(i-a)\epsilon$ for $a\in A$ as above.  
We could view $\epsilon$ as the symbol for an infinitesimal number. 
The weights \eqref{eq:Ellweights} specialize to 
\begin{equation}\label{eq:Ellparaweight}
\begin{matrix}
&\begin{tikzpicture}
\redstrings
\draw[red,thick,<-] (0,0) node [left] {$a$}
    \b\u\e;
\draw[red,thick,<-] (0,1) node [left] {$b$}
    \b\d\e;
\bluestrings
\draw[blue,thick,<-] (0,0.1) 
    \b\u\e node [right] {$d+t\epsilon$};
\draw[blue,thick,<-] (0,1.1) 
    \b\d\e node [right] {$c+s\epsilon$};
\end{tikzpicture}& 
\begin{tikzpicture}
\redstrings
\draw[red,thick,<-] (0,0) node [left] {$a$}
    \b\u\e;
\draw[red,thick,<-] (0,1) node [left] {$b$}
    \b\d\e;
\bluestrings
\draw[blue,thick,<-] (0,0.1) 
    \b\h\e node [right] {$c+s\epsilon$};
\draw[blue,thick,<-] (0,1.1) 
    \b\h\e node [right] {$d+t\epsilon$};
\end{tikzpicture}\\
[\operatorname{weight}]_P&
\Q(\la_{c}\-\la_d\+(t\-s)\hbar,z_a\-z_b)&
\P(\la_{c}\-\la_d\+(t\-s)\hbar,z_a\-z_b).
\end{matrix}
\end{equation}
%
%
%
We can now translate Theorem \ref{th:Billeypara} as follows. 

\begin{Th}\label{th:BilleyparatypeA}
Let $u,w\in W$. 
For a fixed wiring diagram $D_0$ of $u$, we have
$$\bfE^P_w(u)=\sum_{D}[\operatorname{weight}(D)]_P, $$
where the sum is over all sub-wiring diagrams $D$ of $D_0$ whose underlying permutation is $w$. 
We can further assume that any pair of blue strings in $D$ labeled by $c+\mathbb{N}\epsilon$ for the same $c\in A$ does not intersect.
\end{Th}



\section{Pipe Dream Model}
\label{ref:pipedream}

In this section, we consider polynomial representatives of equivariant elliptic Schubert classes $\bfE_w$ in type $A$ via the generic pipe dreams of Knutson and Zinn-Justin \cite{KZ}, which were recalled in Section~\ref{intro}. 

We remark that some polynomial representatives for the elliptic Schubert classes of \cite{RW1,KRW} in type $A$ are obtained by certain weight functions, and we do not know a precise relationship to our model, since our classes are different from theirs; see Theorem \ref{thm:dualbasis}.

\subsection{Polynomial representatives}\label{subs:polrep}
Recall the geometric meaning of $\mathcal{Q}_W^*$ as the localization model of 
\eqref{eq:defEGB}. 
We can also consider the Borel model of $K_T(G/B)$. When $G=\operatorname{GL}_n$, we  consider 
$$\mathcal{P} = \mathbb{Q}[\Lambda][[q]],\qquad 
\Lambda = 
\left(\bigoplus_{i=1}^n\mathbb{Z}x_i\right)\oplus 
\left(\bigoplus_{i=1}^n\mathbb{Z}y_i\right)\oplus 
\left(\bigoplus_{i=1}^n\mathbb{Z}\la_i\right)\oplus 
\mathbb{Z}\hbar.$$
Denote by $x_i,y_i,\la_i$ for $1\leq i\leq n$ the coordinates of the three summands. 
Define the localization map, which is a ring homomorphism:
$$\operatorname{Loc}:\mathcal{P}\longrightarrow \mathcal{Q}_W^*,\qquad 
\operatorname{Loc}(f(x,y,\la))(u) = f(uz,z,\la),$$
where $uz=(z_{u(1)},\ldots,z_{u(n)})$. 
A \emph{polynomial representative} of $\bfE_w$ is an element $\mathcal{E}_w\in \mathcal{P}$ such that $\operatorname{Loc}(\mathcal{E}_w)=\bfE_w$. 
Note that the polynomial representative is not unique. 

\begin{Rmk} {\rm The above notation agrees up to sign with the traditional one for the variables in double Schubert and Grothendieck polynomials, in the sense that $x_i$ are the Chern roots of the universal subquotient bundles on $\operatorname{GL}_n/B$, and $y_i$ are the equivariant parameters (characters of the torus); in addition, we now also have the dynamical variables $\la_i$. More precisely, the pipe dream formula for double Schubert polynomials (see, e.g., \cite[Section 3]{knu}) yields, via specialization, a formula for the localizations of Schubert classes which is positive in $x_j-x_i$, $j>i$ (i.e., the negative roots), while our corresponding formula is positive in the positive roots. }
\end{Rmk}

We construct the  polynomial representatives $\mathcal{E}_w$ for $w$ in $W=S_n$.
The construction is based on the Billey-type formula in a larger group $\widehat{G}=\operatorname{GL}_{2n}$, whose Weyl group is $\widehat{W}= S_{2n}$. 
 Consider $\widehat{W}_{\widehat{P}}=S_1\times \cdots\times S_1\times S_n$, and let 
$u_0\in S_{2n}$ be given by 
$$1\leq i\leq n,\qquad u_0(i)=n+i,\qquad u_0(n+i)=i.$$
For $w\in W=S_n$, we  view $w=w\times \id\in S_n\times S_n\subset  S_{2n}$. 
From the discussion in the previous sections, 
we have constructed 
$$
\bfE_{w}^{\hat{P}}(u_0)=[\bfE_{w}(u_0)]_{\widehat{P}}\in 
\mathbb{Q}[\Lambda][[q]],\qquad 
\Lambda = 
\left(\bigoplus_{i=1}^{2n}\mathbb{Z}z_i\right)\oplus 
\left(\bigoplus_{i=1}^{n+1}\mathbb{Z}\la_i\right)\oplus 
\mathbb{Z}\hbar.$$
In particular, we can write it as 
$$\bfE_{w}^{\widehat{P}}(u_0)=
\widehat{\mathcal{E}}_{w}(z_1,\ldots,z_{2n};\la_{1},\ldots,\la_{n},\la_{n+1}).$$
We  define 
$$\mathcal{E}_w(x,y,\la)
=\widehat{\mathcal{E}}_{w}(y_1,\ldots,y_n,x_1,\ldots,x_n;\la_{1},\ldots,\la_{n},0)\in \mathcal{P}.$$

\begin{Lemma}\label{lem:LrecforPolyE}
For a simple root $\al$ of $\operatorname{GL}_n$, we have
\begin{equation}\label{eq:LrecforPolyE}
\mathcal{E}_{s_\al w}(x,y,\la)
= 
-\frac{\P(\la_{w^{\-1}\al^\v},\-y_\al)}
    {\Q(\-\la_{w^{\-1}\al^\v},\-y_\al)}
    \mathcal{E}_{w}(x,y,\la)+
\frac{1}
    {\Q(\-\la_{w^{\-1}\al^\v},\-y_\al)}
    \mathcal{E}_{w}(x,s_\al y,\la).
\end{equation}
\end{Lemma}
\begin{proof}
Note that $s_\al uW_P=uW_P$. 
Thus by \eqref{lem:buw_para2} of Lemma \ref{lem:buw_para} and \eqref{eq:rec-bleft} in Proposition \ref{prop:rec-b}, we have 
$$\mathcal{E}_{w}(z,\la)= 
    \P(\la_{w^{\-1}\alv},z_{\al})\cdot
    \mathcal{E}_{w}(s_\al z,\la) +
    \Q(\-\la_{w^{\-1}\alv},z_{\al})\cdot 
    \mathcal{E}_{s_\al w}(s_\al z,\la).
$$
By specializing $z_1,\ldots,z_{2n}$ to $y_1,\ldots,y_n,x_1,\ldots,x_n$, the identity \eqref{eq:LrecforPolyE} follows. 
\end{proof}

Let $w_\circ$ be the longest element of $W=S_n$.

\begin{Lemma}\label{lem:Ew0init}
We have 
\begin{equation}\label{eq:Ew0init}
\mathcal{E}_{w_\circ}(x,y,\la)=\prod_{i+j\leq n}
\Q(\la_i-\la_{n+1-j},y_{j}-x_i)
\prod_{i=1}^n\P(\la_i,y_{n+1-i}-x_{i}).
\end{equation}
\end{Lemma}
\begin{proof}
This is immediate based on the translation from sub-wiring diagrams to generic pipe dreams, which is explained in Section~\ref{subsec:genpd} below. 
More precisely, there is only one generic pipe dream for $w_\circ$, and its weight is given by  \eqref{eq:Ew0init}. 
The case $n=4$ is illustrated below.
$$\BPD{
\M{1}\X\X\X\J[0]\\
\M{2}\X\X\J[0]\O\\
\M{3}\X\J[0]\O\O\\
\M{4}\J[0]\O\O\O}
\qquad \begin{matrix}
\Q(\la_1\-\la_4,y_1\-x_1)&
    \Q(\la_1\-\la_3,y_2\-x_1)&
        \Q(\la_1\-\la_2,y_3\-x_1)&  
            \P(\la_1,y_4\-x_1)\\[1ex]
\Q(\la_2\-\la_4,y_1\-x_2)&
    \Q(\la_2\-\la_3,y_2\-x_2)&
        \P(\la_2,y_3\-x_2)&
            1\\[1ex]
\Q(\la_3\-\la_4,y_1\-x_3)&
    \P(\la_3,y_2\-x_3)&
        1&
            1\\[1ex]
\P(\la_4,y_1\-x_4)&
    1&
        1&
            1\\
\end{matrix}$$
\end{proof}

\begin{Th}
The element $\mathcal{E}_w\in \mathcal{P}$ is a polynomial representative of $\bfE_w$. 
\end{Th}
\begin{proof}
We first check 
$$\bfE'_{w_\circ}:=\operatorname{Loc}(\mathcal{E}_{w_\circ}(x,y,\la))=\bfE_{w_\circ}$$
using Lemma \ref{lem:Ew0init}. 
By \eqref{eq:PQnorm}, we have 
$$\bfE'_{w_\circ}(u)=
\prod_{i+j\leq n}
\Q(\la_i-\la_{n+1-j},z_{j}-z_{u(i)})
\prod_{i=1}^n\P(\la_i,z_{n+1-i}-z_{u(i)})=0$$
unless $u(i)\neq j$ for any $j\leq n-i$, i.e. $u=w_\circ$. 
For $u=w_\circ$, we have 
$$\bfE'_{w_\circ}(w_\circ)=
\prod_{i+j\leq n}
\Q(\la_i-\la_{n+1-j},z_{j}-z_{n+1-i})
$$
by \eqref{eq:PQnorm}. This agrees with the computation of $\bfE_{w_\circ}(w_\circ)$ in Example \ref{eg:auubuu}.

In general, Lemma \ref{lem:LrecforPolyE} and \eqref{eq:rec-bleft} in Proposition \ref{prop:rec-b} implies that 
$$
\operatorname{Loc}(\mathcal{E}_w)=\bfE_w
\Longrightarrow 
\operatorname{Loc}(\mathcal{E}_{s_\al w})=\bfE_{s_\al w}.
$$
As a result, the theorem follows by downward induction. 
\end{proof}

\subsection{Generic pipe dreams}\label{subsec:genpd}
From our definition and Theorem \ref{th:BilleyparatypeA}, there is a sub-wiring diagram model for $\mathcal{E}_w$. 
Note that the wiring diagrams for any reduced word of $u_0$ are isotopically equivalent. 
The set of sub-wiring diagrams is in a natural bijection with the set of generic pipe dreams \cite{KZ} by replacing the {\color{red}red} strings by the ``Poincar\'e dual'', i.e., an $n\times n$ grid,  and the {\color{blue}blue} strings labeled by $1,\ldots,n$ by pipes in the grid. 
Since the {\color{blue}blue} strings labeled by $n+1+\mathbb{N}\epsilon$ are not allowed to intersects, by Theorem \ref{th:BilleyparatypeA}, we lose no information. 

Here is an example of this bijection for $n=3$.
$$
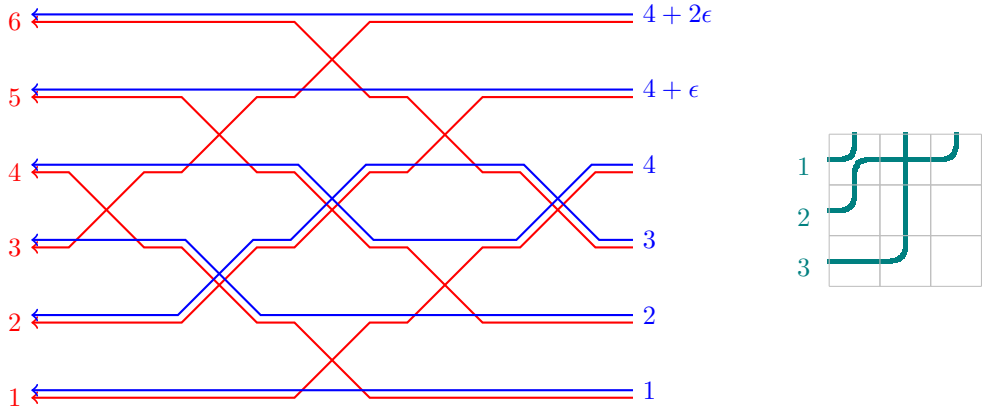
\begin{matrix}
\begin{tikzpicture}
\redstrings
\draw[red,thick,<-](0,1) node [left] {1} 
    \b\h\h\u\u\u\e;
\draw[red,thick,<-](0,2) node [left] {2} 
    \b\h\u\u\u\h\e;
\draw[red,thick,<-](0,3) node [left] {3} 
    \b\u\u\u\h\h\e;
\draw[red,thick,<-](0,4) node [left] {4} 
    \b\d\d\d\h\h\e;
\draw[red,thick,<-](0,5) node [left] {5} 
    \b\h\d\d\d\h\e;
\draw[red,thick,<-](0,6) node [left] {6} 
    \b\h\h\d\d\d\e;
\bluestrings
\draw[blue,thick,<-] (0,1.1)
    \b\h\h\h\h\h\e node [right] {1};
\draw[blue,thick,<-] (0,2.1) 
    \b\h\u\u\h\d\e node [right] {3};
\draw[blue,thick,<-] (0,3.1) 
    \b\h\d\h\h\h\e node [right] {2};
\draw[blue,thick,<-] (0,4.1) 
    \b\h\h\d\h\u\e node [right] {4};
\draw[blue,thick,<-] (0,5.1)
    \b\h\h\h\h\h\e node [right] {$4+\epsilon$};
\draw[blue,thick,<-] (0,6.1)
    \b\h\h\h\h\h\e node [right] {$4+2\epsilon$};
\end{tikzpicture}
\end{matrix}
\qquad 
\BPD[1.6pc]{
\M{1}\B\X\J\\\M{2}\J\I\O\\\M{3}\H\J\O
}$$

We now explain the translation of the weights to generic pipe dreams. We first need to define the \emph{level} statistic.

\begin{Def}\label{def:level}
For a generic pipe dream, we define \emph{level} for tiles containing only one pipe inside, i.e. $\BPD{\F},\BPD{\J},\BPD{\I}, \BPD{\H}$, as follows. 
For each pipe, imagine that we are walking along  it with a number kept in mind. At the source of the pipe, the number is $0$. 
Then once we meet 
$\BPD{\F},\BPD{\J},\BPD{\I}, \BPD{\H}$
we define its \emph{level} to be the number in mind, and then 
increase the number by $1$ if we met $\BPD{\H}$,
and by $-1$ if met $\BPD{\I}$. 
\end{Def}

We can now translate the weights at the $(i,j)$-position as follows.
$$
\begin{array}{cc}
\BPD[1.2pc]{\X} & \Q(\la_a\-\la_b,y_j\-x_i) \\[2ex]
\BPD[1.2pc]{\B} & \P(\la_a\-\la_b,y_j\-x_i)
\end{array}\quad
\begin{array}{cc}
\BPD[1.2pc]{\H} & \Q(\la_a\-c\hbar,y_i\-x_i)
\\[2ex]
\BPD[1.2pc]{\J} & \P(\la_a\-c\hbar,y_i\-x_i)
\end{array}\quad
\begin{array}{cc}
\BPD[1.2pc]{\I} & \Q(c\hbar\-\la_b,y_i\-x_i)
\\[2ex]
\BPD[1.2pc]{\F} & \P(c\hbar\-\la_b,y_i\-x_i)
\end{array}
\quad 
\begin{array}{cc}
\BPD[1.2pc]{\O} & 1 
\end{array}$$
Here $a$ (resp., $b$) is the index for the pipe from the lower (resp., left) boundary, and $c$ is the level. 
See the following example, where the levels are displayed inside the respective tiles. 
$$
\operatorname{weight}\left(\BPD{
\M{1}\B\X\J[0]\\\M{2}\J[0]\I[0]\O\\\M{3}\H[0]\J[1]\O
}\right)=
\prod\left\{
\begin{matrix}
\P(\la_1\-\la_2,y_1\-x_1) & 
    \Q(\la_2\-\la_3,y_2\-x_1) & 
        \P(\la_2,y_3\-x_1)\\[1ex]
\P(\la_2,y_1\-x_2) & 
    \Q(\-\la_3,y_2\-x_2) & 
        1\\[1ex]
\Q(\la_3,y_1\-x_3) & 
    \P(\la_3\-\hbar,y_2\-x_3) & 
        1\\
\end{matrix}\right\}$$

\begin{Th}\label{thm:pipedreamofEw}
For any $w\in S_n$, we have 
$$\mathcal{E}_w = \sum_{\pi\in \mathsf{GPD}(w)}\operatorname{weight}(\pi).$$
\end{Th}

\begin{Rmk}{\rm We can define polynomial representatives of equivariant elliptic classes for partial flag manifolds $\operatorname{GL}_n/P$ by applying the map $[\_]_P$ to the polynomial representatives for $\operatorname{GL}_n/B$. In this way, we derive the analogue of Theorem~\ref{thm:pipedreamofEw} for $\operatorname{GL}_n/P$; the respective sum is now restricted to those generalized pipe dreams for which the strings labeled by $c+\mathbb{N}\epsilon$ for the same $c\in A$ do not intersect (see the notation in Theorem~\ref{th:BilleyparatypeA}), and the map $[\_]_P$ is applied to the corresponding weights (as in Theorem~\ref{th:BilleyparatypeA}). 
}
\end{Rmk}


\appendix
\section{Limit to $K$-theory}
\label{sec:Kthylimit}

In this appendix, we study the degeneration of elliptic Schubert classes to $K$-theory. A similar computation also appeared in \cite[Proposition 12]{Smi20} (see also \cite{KS22, KS23}). 


\begin{Prop}\label{Prop:limtheta}Denote $e(u)=e^{2\pi iu}$. 
Let $\tau$ be a variable such that $e(\tau)=q$. We have 
\begin{align}
\lim_{q\to 0} \theta(x) & = e(\tfrac{x}{2})-e(\-\tfrac{x}{2}), &\lim_{q\to 0} \frac{\theta(u+s\tau)}{\theta(s\tau)} & = e\big((\lfloor -s\rfloor+\tfrac{1}{2}) u\big), \quad s\in \mathbb{Q}\setminus \mathbb{Z}. 
\end{align}
\end{Prop}

\begin{proof}
The first limit is well-known. For the second one, we consider the case  $s>0$ since the case $s<0$ is similar.
Denote $x=e(u)=e^{2\pi iu}$. Then we have 
\begin{align*}
\frac{\theta(u+s\tau)}{\theta(s\tau)}
& = 
\frac
{q^{s/2} x^{1/2}-q^{-s/2}x^{-1/2}}
{q^{s/2}-q^{-s/2}}
\prod_{n>0}
\frac
{(1-q^{n+s}x)(1-q^{n-s}/x)}
{(1-q^{n+s})(1-q^{n-s})}\\
& = 
\frac
{q^{s/2} x^{1/2}-q^{-s/2}x^{-1/2}}
{q^{s/2}-q^{-s/2}}
\prod_{n>0}
\frac
{1+q^{2n}-q^{n+s}x-q^{n-s}/x}
{1+q^{2n}-q^{n+s}-q^{n-s}}.
\end{align*}
It suffices to compute the limit of each factor. 
When $n>s$, we have $2n>n+s>n-s>0$, so 
$$\lim_{q\to 0}
\frac
{1+q^{2n}-q^{n+s}x-q^{n-s}/x}
{1+q^{2n}-q^{n+s}-q^{n-s}}=1.$$
When $n<s$, we have $n+s>2n>0>n-s$, so 
$$\lim_{q\to 0}\frac
{1+q^{2n}-q^{n+s}x-q^{n-s}/x}
{1+q^{2n}-q^{n+s}-q^{n-s}}=x^{-1}.$$
Since 
$$\lim_{q\to 0}
\frac
{q^{s/2} x^{1/2}-q^{-s/2}x^{-1/2}}
{q^{s/2}-q^{-s/2}}
=x^{-1/2}.$$
So the total contribution is $e\big((-\lfloor s\rfloor-\tfrac{1}{2}) u\big)=e\big((\lfloor-s\rfloor+\tfrac{1}{2}) u\big)$.
\end{proof}

\begin{Rmk}{\rm 
In our paper, we treat $\theta$-functions formally \eqref{eq:convtheta}, so we need a few words explaining the meaning of the limit. 
For a lattice $\Lambda$, we have 
$
\operatorname{Frac}\left(\mathbb{Q}[\Lambda][[q]]\right)\subseteq
\left(\operatorname{Frac}\mathbb{Q}[\Lambda]\right)(\!(q)\!)$. 
So the meaning of the limit is obvious for $f\in \operatorname{Frac}\left(\mathbb{Q}[\Lambda][[q]]\right)$, i.e., we have }
$$\lim_{q\to 0}f(q) = \begin{cases}
f(0), & 
f\in \left(\operatorname{Frac}\mathbb{Q}[\Lambda]\right)[[q]],\\
\mathsf{DNE}, & \mbox{otherwise}. 
\end{cases}$$
{\rm More generally, functions $f$ such as $\theta(u+s\tau)$
for some $s\in \mathbb{Q},u\in \Lambda$ take values  in the Puiseux field 
$
\bigcup_{N>0}\left(\operatorname{Frac}\mathbb{Q}[\Lambda]\right)(\!(q^{1/N})\!)$.
Then for such an element $f$, we have}
$$\lim_{q\to 0}f(q) = \begin{cases}
f(0), & 
f\in \bigcup_{N>0}\left(\operatorname{Frac}\mathbb{Q}[\Lambda]\right)[[q^{1/N}]],\\
\mathsf{DNE}, & \mbox{otherwise}. 
\end{cases}$$
\end{Rmk}

\begin{Coro}\label{coro:limitPQ}
For $0<s<1$, taking the limit as $q\to 0$, we have 
\begin{align}
\P(s\tau,y)=\four{s\tau-y}{\hbar}{y+\hbar}{s\tau}
& \to \frac{e(\tfrac{y}{2})\big(e(\tfrac{\hbar}{2})-e(\-\tfrac{\hbar}{2}))\big)}{e(\tfrac{y+\hbar}{2})-e(\-\tfrac{y+\hbar}{2}))}
= \frac{1-e(\-\hbar)}{1-e(\-y)e(\-\hbar)}, \\
\Q(s\tau,y)=\four{s\tau+\hbar}{y}{y+\hbar}{s\tau}&
\to \frac{e(\-\tfrac{\hbar}{2})\big(e(\tfrac{y}{2})-e(\-\tfrac{y}{2}))\big)}{e(\tfrac{y+\hbar}{2})-e(\-\tfrac{y+\hbar}{2}))}
= \frac{(1-e(\-y))e(\-\hbar)}{1-e(\-y)e(\-\hbar)},\\
\P(\-s\tau,y)=\four{-s\tau-y}{\hbar}{y+\hbar}{-s\tau}
& \to \frac{e(\-\tfrac{y}{2})\big(e(\tfrac{\hbar}{2})-e(\-\tfrac{\hbar}{2}))\big)}{e(\tfrac{y+\hbar}{2})-e(\-\tfrac{y+\hbar}{2}))}
= \frac{(1-e(\-\hbar))e(\-y)}{1-e(\-y)e(\-\hbar)},\\
\Q(-s\tau,y)=\four{-s\tau+\hbar}{y}{y+\hbar}{-s\tau}&
\to \frac{e(\tfrac{\hbar}{2})\big(e(\tfrac{y}{2})-e(\-\tfrac{y}{2}))\big)}{e(\tfrac{y+\hbar}{2})-e(\-\tfrac{y+\hbar}{2}))}
= \frac{1-e(\-y)}{1-e(\-y)e(\-\hbar)}.
\end{align}
\end{Coro}


\begin{Def}
Denote $y=-e(-\hbar)$. 
For $0<s\ll 1$ and $w\in W$, we define 
$$\mathcal{Q}^K=\operatorname{Frac}\left(\mathbb{Q}[X^*(T)]\right)(y),\qquad
K_w\in \mathcal{Q}^{K,*}_W=\operatorname{Map}(W,\mathcal{Q}^K)$$
by 
$$K_w(u)=y^{-\ell(w)}\lim_{q\to 0} (\bfE_w(u)\big|_{\displaystyle
\begin{subarray}{l}
\la_\alv=s\tau\\
\forall\, \al\in \Sigma
\end{subarray}})
\in \mathcal{Q}^{K,*}_W.$$
\end{Def}

\begin{Rmk}\label{rmk:EKisSMC}{\rm 
If $s_\al w>w$, i.e., $w^{-1}\al>0$, then $\la_{w^{\-1}\al^\v}$ specializes to $hs$, with  $h>0$ being the height of $w^{\-1}\alpha$. 
By \eqref{eq:rec-bleft} of Proposition \ref{prop:rec-b} and the limit in Corollary \ref{coro:limitPQ}, we have }
$$K_w(s_\al u) = 
    \frac{1+y}{1+ye(-\al)}\cdot
    {}^{s_\al}K_w(u)+
    (-y)\frac{1-e(-\al)}{1+ye(-\al)}\cdot
    {}^{s_\al}K_{s_\al w}(u).
$$
{\rm Similarly, if  $s_\al w<w$, i.e., $w^{-1}\al<0$, we have}
$$K_w(s_\al u) = 
    \frac{(1+y)e(-\al)}{1+ye(-\al)}\cdot
    {}^{s_\al}K_w(u)+
    (-y)^{-1}\frac{-y(1-e(-\al))}{1+ye(-\al)}\cdot
    {}^{s_\al}K_{s_\al w}(u).
$$
{\rm From this, we can conclude that
$
K_w=\mathsf{SMC}_{y}(X^{\circ,w})$; see \cite[Theorem 7.8]{MNS}.}
\end{Rmk}

Similarly, for $w\in W^P$, we can define the parabolic version 
$$K_w^{P}(u)
=y^{-\ell(w)}\lim_{q\to 0}
\bfE_w^P(u)\big|_{\displaystyle
\begin{subarray}{l}
\la_\alv=s\tau\\
\forall\, \al\notin \Sigma_P
\end{subarray}}
=y^{-\ell(w)}\lim_{q\to 0}
\bfE_w(u)\big|_{\displaystyle
\begin{subarray}{l}
\la_\alv=-\hbar ,\forall \al\in \Sigma_P\\
\la_\alv=s\tau ,\forall \al\notin \Sigma_P
\end{subarray}}\,,
$$
and it is not hard to show that
$K^{P}_w=\mathsf{SMC}_{y}(X^{\circ,wW_P})$.

We now describe the limit in type $A$. 
\begin{itemize}
    \item 
The specialization $\la_\alv=s\tau$ is equivalent to the specialization $\la_i=\la_1-s\cdot i$. 
As a result, the limit splits into two cases by comparing the input blue strings. 

    \item 
The normalization $y^{-\ell(w)}$ can be absorbed into weights via multiplication by $(-y)^{\pm 1}$, depending on the cross being positive or negative. 
\end{itemize}
The resulting weights obtained from  \eqref{eq:Ellparaweight} 
are listed below.
$$\begin{matrix}
&\begin{tikzpicture}
\redstrings
\draw[red,thick,<-] (0,0) node [left] {$a$}
    \b\u\e;
\draw[red,thick,<-] (0,1) node [left] {$b$}
    \b\d\e;
\bluestrings
\draw[blue,thick,<-] (0,0.1) 
    \b\u\e node [right] {$d+t\epsilon$};
\draw[blue,thick,<-] (0,1.1) 
    \b\d\e node [right] {$c+s\epsilon$};
\end{tikzpicture}& 
\begin{tikzpicture}
\redstrings
\draw[red,thick,<-] (0,0) node [left] {$a$}
    \b\u\e;
\draw[red,thick,<-] (0,1) node [left] {$b$}
    \b\d\e;
\bluestrings
\draw[blue,thick,<-] (0,0.1) 
    \b\h\e node [right] {$c+s\epsilon$};
\draw[blue,thick,<-] (0,1.1) 
    \b\h\e node [right] {$d+t\epsilon$};
\end{tikzpicture}\\
\operatorname{weight}&
\begin{array}{cl}
(-y)^{-1}\dfrac{-y(1-e(-\al))}{1+ye(-\al)}, &  c<d\\[2ex]
(-y)^1\dfrac{1-e(-\al)}{1+ye(-\al)}, &  c>d\\[2ex]
1, &  c=d
\end{array}&
\begin{array}{rl}
\dfrac{1+y}{1+ye(-\al)}, &  c<d\\[2ex]
\dfrac{(1+y)e(-\al)}{1+ye(-\al)}, &  c>d\\[2ex]
0, &  c=d
\end{array}
\end{matrix}
$$
If we further set $y=0$, Theorem \ref{thm:pipedreamofEw} specializes to the pipe dream model for double Grothendieck polynomials \cite{FK}. 


\section{Comparing Elliptic Classes}
\label{sec:comparing}

In this section, we discuss the comparison with the classes introduced in \cite{RW1} and \cite{KRW}. 
Recall the Poincar\'e pairing \eqref{eq:PoincarePair}. 
In \cite[Section 1.8]{RW1}, for the Schubert variety $X_w$, the \emph{localized elliptic class} $E(X_w)$ are defined whose restrictions are denoted by $E(X_w)_{u}\in \calQ$, and  many examples are computed in \cite[Section 12]{RW1}.   
We define the \emph{rescaled elliptic classes} $\mathfrak{E}(X_w),\mathfrak{E}'(X_w)\in \mathcal{Q}_W^*$ for $w\in W$ by 
\begin{align}
\label{eq:rescalefrakEXw}
\mathfrak{E}(X_w)(u) &= \prod_{\begin{subarray}{c}\al>0\\w\al<0\end{subarray}}\frac{\theta(\la_\alv)\theta(\hbar)}{\theta(\la_\alv-\hbar)\theta'(0)}\cdot\prod_{\alpha>0}\theta(-z_{u\al})\cdot
 E(X_w)_u,~u\le w,\\
\label{eq:renormfrakEXw}
\mathfrak{E}'(X_w)(u) &= \prod_{\begin{subarray}{c}\al>0\\w\al<0\end{subarray}}\frac{\theta(\la_\alv)\theta(\hbar)}{\theta(\la_\alv+\hbar)\theta'(0)}
 \cdot\prod_{\alpha>0}\theta(-z_{u\al})\cdot E(X_w)_u, ~u\le w. 
\end{align}
Note that in \cite{RW1},  \eqref{eq:renormfrakEXw} is the definition, but the $R$-matrix recursion is actually for \eqref{eq:rescalefrakEXw}. Indeed, 
the classes $\mathfrak{E}(X_w)$ are characterized by the initial condition
$$ \mathfrak{E}(X_\id) = \left({\prod_{\al>0}\theta(-z_\al)}\right)f_{\id}$$
and the ``$R$-matrix recursion'' \cite[Theorem 1.3]{RW1}:
\begin{align*}
\mathfrak{E}(X_{s_\al w})(u)
& =\frac{\delta(-z_{\al},\la_{w^{\-1}\alv})}
        {\delta(-\la_{w^{\-1}\alv},\hbar)}
    \mathfrak{E}(X_w)(u)+
    \frac{\delta(z_\al,\hbar)}
        {\delta(-\la_{w^{\-1}\alv},\hbar)}
    {}^{s_\al}\mathfrak{E}(X_w)(s_\al u)\\
& = 
-\frac{\P(\la_{w^{\-1}\al^\v},z_\al)}
    {\Q(\-\la_{w^{\-1}\al^\v},z_\al)}
    \mathfrak{E}(X_w)(u)+
\frac{1}
    {\Q(\-\la_{w^{\-1}\al^\v},z_\al)}
    {}^{s_\al}\mathfrak{E}(X_w)(s_\al u).
\end{align*}

For a class $\gamma\in \calQ_{W}^*$ and $w\in W$, we denote by ${}^w\gamma\in \calQ_W^*$ the class with $({}^w\gamma)(u)={}^w(\gamma(w^{-1}u))$.
It is easy to check that
\begin{equation}\label{eq:1}
\langle {}^w \gamma_1,{}^w\gamma_2\rangle
={}^w \langle \gamma_1,\gamma_2\rangle. 
\end{equation} 
Using this notation, the above recursion can be written as 
\begin{equation}\label{eq:LrecforEX}
\mathfrak{E}(X_{s_\al w})
= 
-\frac{\P(\la_{w^{\-1}\al^\v},z_\al)}
    {\Q(\-\la_{w^{\-1}\al^\v},z_\al)}
    \mathfrak{E}(X_w)+
\frac{1}
    {\Q(\-\la_{w^{\-1}\al^\v},z_\al)}
    {}^{s_\al}\mathfrak{E}(X_w).
\end{equation}
Similarly, \eqref{eq:rec-bleft} can be written as 
\begin{equation}\label{eq:LrecforE}
\bfE_{s_\al w}
 = 
-\frac{\P(\la_{w^{\-1}\al^\v},\-z_\al)}
    {\Q(\-\la_{w^{\-1}\al^\v},\-z_\al)}
    \bfE_{w}+
\frac{1}
    {\Q(\-\la_{w^{\-1}\al^\v},\-z_\al)}
    {}^{s_\al}\bfE_{w}.
\end{equation}

\begin{Th}\label{thm:dualbasis}
For any $u,w\in W$, 
we have 
\begin{equation}\label{eq:EXuEwpair}
\langle\mathfrak{E}(X_u),\bfE_w\rangle
=\delta_{u,w}^{\mathsf{Kr}}\prod_{
\begin{subarray}{c}
\al>0\\
u\al<0
\end{subarray}}
\frac{\theta(\la_\alv+\hbar)}{\theta(\la_\alv-\hbar)},
\qquad 
\langle\mathfrak{E}'(X_u),\bfE_w\rangle
=\delta_{u,w}^{\mathsf{Kr}}. 
\end{equation}
In other words, the elliptic Schubert classes
$\{\bfE_w\}$ form the dual basis to $\{\mathfrak{E}'(X_w)\}$ via the Poincar\'e pairing \eqref{eq:PoincarePair}. 
\end{Th}
\begin{proof}
It is obvious that the second equation follows from the first, so we prove the first equation by induction on $u$. 
If $u=\id$, by definition, we have 
$$\langle \mathfrak{E}(X_{\id}),\gamma\rangle = \gamma(\id),\qquad 
\text{for any }\gamma\in \calQ_{W}^*.$$
Since $\bfE_w(\id)=\delta_{w,\id}^{\mathsf{Kr}}b_{\id, \id}$, \eqref{eq:EXuEwpair} follows for $u=\id$. 

Assuming that \eqref{eq:EXuEwpair} is true for $u$,  that is, the pairing of $\frakE(X_u)$ with $\bfE_w$ satisfies the equations for all $w\in W$. We will  compute $\langle \frakE(X_{s_\al u}), \bfE_w\rangle$.  
By pairing \eqref{eq:LrecforE} with $\mathfrak{E}(X_u)$, we get 
\begin{align*}
\langle \mathfrak{E}(X_u),\bfE_{s_\al w}\rangle
 = 
-\frac{\P(\la_{w^{\-1}\al^\v},\-z_\al)}
    {\Q(\-\la_{w^{\-1}\al^\v},\-z_\al)}
\langle \mathfrak{E}(X_u),\bfE_{w}\rangle +
\frac{1}
    {\Q(\-\la_{w^{\-1}\al^\v},\-z_\al)}
\langle \mathfrak{E}(X_u),{}^{s_\al}\bfE_{w}\rangle,
\end{align*}
so by induction assumption, we have 
\begin{equation}
\label{eq:3}\langle \mathfrak{E}(X_u),{}^{s_\al}\bfE_{w}\rangle
=\left(\begin{cases}
\P(\la_{u^{\-1}\al^\v},\-z_\al),& w=u,\\
\Q(\la_{u^{\-1}\al^\v},\-z_\al),& w=s_\al u,\\
0, & \text{otherwise}
\end{cases}\right)\cdot
\langle \mathfrak{E}(X_{u}),\bfE_{u}\rangle.
\end{equation}
In particular, 
$\langle \mathfrak{E}(X_u),{}^{s_\al}\bfE_{w}\rangle= 0$ unless $w\in \{u,s_\al u\}$. 
By pairing \eqref{eq:LrecforEX} with $\bfE_w$, we get 
\begin{align*}
\langle\mathfrak{E}(X_{s_\al u}),\bfE_{w}\rangle 
& = 
-\frac{\P(\la_{u^{\-1}\al^\v},z_\al)}
    {\Q(\-\la_{u^{\-1}\al^\v},z_\al)}
    \langle \mathfrak{E}(X_u),\bfE_{w}\rangle+
\frac{1}
    {\Q(\-\la_{u^{\-1}\al^\v},z_\al)}
    \langle {}^{s_\al}\mathfrak{E}(X_u),\bfE_{w}\rangle\\
& \overset{\eqref{eq:1}}=
-\frac{\P(\la_{u^{\-1}\al^\v},z_\al)}
    {\Q(\-\la_{u^{\-1}\al^\v},z_\al)}
    \langle \mathfrak{E}(X_u),\bfE_{w}\rangle+
\frac{1}
    {\Q(\-\la_{u^{\-1}\al^\v},z_\al)}
    {}^{s_\al}\langle \mathfrak{E}(X_u), {}^{s_\al} \bfE_{w}\rangle.
\end{align*}
Thus by \eqref{eq:3} and the induction assumption, $\langle\mathfrak{E}(X_{s_\al u}),\bfE_{w}\rangle=0$ unless $w\in \{u,s_\al u\}$. 
If $w=u$, then we have
$$
\langle\mathfrak{E}(X_{s_\al u}),\bfE_{w}\rangle \overset{\eqref{eq:3}}=
-\frac{\P(\la_{u^{\-1}\al^\v},z_\al)}
    {\Q(\-\la_{u^{\-1}\al^\v},z_\al)}
    \langle \mathfrak{E}(X_u),\bfE_{u}\rangle
+
\frac{\P(\la_{u^{\-1}\al^\v},z_\al)}
    {\Q(\-\la_{u^{\-1}\al^\v},z_\al)}
{}^{s_\al}    \langle \mathfrak{E}(X_u),\bfE_{u}\rangle=0.
$$
Here the last identity we also used the fact that $\langle\frakE(X_u), \bfE_u\rangle$ is invariant under the action of $s_\al$. 
If $w=s_\al u$, then we have
$$
\langle\mathfrak{E}(X_{s_\al u}),\bfE_{w}\rangle 
\overset{\eqref{eq:3}}=
\frac{\Q(\la_{u^{\-1}\al^\v},z_\al)}
    {\Q(\-\la_{u^{\-1}\al^\v},z_\al)}
\langle \mathfrak{E}(X_{u}),\bfE_{u}\rangle
=\frac{\theta(\la_{u^{\-1}\alv}+\hbar)}
{\theta(\la_{u^{\-1}\alv}-\hbar)} 
\langle \mathfrak{E}(X_{u}),\bfE_{u}\rangle.$$
The proof is finished. 
\end{proof}
Consequently, we have
\begin{Coro} For any $u,w\in W$, we have
\begin{equation}\label{eq:6}
\sum_{v\in W} E(X_u)_v\cdot \bfE_w(v)=\prod_{\al>0, u\al<0}\frac{\theta(\la_\alv+\hbar)\theta'(0)}{\theta(\la_\alv)\theta(\hbar)}\de_{w,u}^{\mathsf{Kr}}. 
\end{equation}
\end{Coro}
Next, we can compare with the parabolic elliptic Schubert class studied in \cite{KRW}.
In this case, the elliptic class takes values in a localization of
$$\mathcal{M}^P=K_T(G/P)[H_2(G/P;\mathbb{Z})\oplus\mathbb{Z}\hbar][[q]].$$
By identifying $(G/P)^T=W^P$, $K_T(\pt)=\mathbb{Q}[X^*(T)]$, and $H_2(G/P;\mathbb{Z})=\check{Q}^P$, we can embed $\mathcal{M}^P\subseteq \mathcal{Q}_{W/W_P}^*$ by localization. 
The natural \emph{elliptic Poincar\'e pairing} is defined to be
\begin{equation}\label{eq:paraPoincare}
 \langle f,g\rangle_{G/P} =\sum_{u\in W^P}
\frac{f(u)g(u)}{\prod_{\al\in \Phi^+\setminus \Phi_P^+}\theta(-z_{u\al})}\in \mathcal{Q}
, ~\text{ where }f,g\in \mathcal{Q}_{W/W_P}^*.
\end{equation}
According to \cite[Corollary 6.3]{KRW}, for $w\in W^P$ with corresponding Schubert variety $X_w^P$ for $G/P$, the localized elliptic Schubert classes are related by 
\begin{equation}
\label{eq:7}E(X_w^P)_{u}=\sum_{v\in W_P}[E(X_{w})_{uv}]_P, ~u\in W^P.
\end{equation}

Similar as the complete flag variety case, we define the rescaled elliptic class $\mathfrak{E}'(X_w^P)\in \calQ_{W/W_P}^*$ for $w\in W^P$ as in \eqref{eq:renormfrakEXw}, namely:
\begin{equation}\label{eq:4}
\mathfrak{E}'(X_w^P)(u) = 
\left[\prod_{\begin{subarray}{c}\al>0\\w\al<0\end{subarray}}\frac{\theta(\la_\alv)\theta(\hbar)}{\theta(\la_\alv+\hbar)\theta'(0)}\right]_P\cdot \prod_{\alpha\in \Phi^+\setminus \Phi_P^+}\theta(-z_{u\al})\cdot 
E(X^P_w)_u, ~u\in W^P. 
\end{equation}
Note that $w\in W^P$ implies that $\{\al>0|w\al<0\}\cap \Sigma_P=\emptyset$, so the denominator in the first term in \eqref{eq:4} does not contain $\theta(\la_\alv+\hbar)$ with $\alv\in \Sigma_P$, so applying the map $[\_]_P$ makes sense. 

\begin{Coro}\label{coro:paradualbasis}
For any $u,w\in W^P$, we have 
\begin{equation}\label{eq:EXuEwpairp}
\langle\mathfrak{E}'(X^P_u),\bfE^P_w\rangle
=\delta_{u,w}^{\mathsf{Kr}}. 
\end{equation}
In other words, the parabolic elliptic Schubert classes 
$\{\bfE^P_w\}$ form the dual basis to $\{\mathfrak{E}'(X^P_w)\}$ via the Poincar\'e pairing \eqref{eq:paraPoincare}. 
\end{Coro}
\begin{proof}
We have 
\begin{align*}
\langle \mathfrak{E}'(X_u^P),\bfE_w^P\rangle_{G/P}
&\overset{\eqref{eq:paraPoincare}}=\sum_{w'\in W^P}\frac{\frakE'(X_u^P)(w')\cdot \bfE_w^P(w')}{\prod_{\al\in \Phi^+\setminus \Phi^+_P}\theta(-z_{w'\al})}\\
&\overset{\eqref{eq:4}}=\sum_{w'\in W^P}\frac{\prod_{\al>0,u\al<0}\frac{\theta(\la_{\alv})\theta(\hbar)}{\theta(\la_\alv+\hbar)\theta'(0)}
\prod_{\al\in \Phi^+\setminus \Phi^+_P}\theta(-z_{w'\al})\cdot E(X^P_u)_{w'}\cdot \bfE_w^P(w')}{\prod_{\al\in \Phi^+\setminus \Phi^+_P}\theta(-z_{w'\al})}\\
&\overset{\eqref{eq:7}}=\sum_{w'\in W^P}\prod_{\al>0,u\al<0}\frac{\theta(\la_{\alv})\theta(\hbar)}{\theta(\la_\alv+\hbar)\theta'(0)}
\sum_{v\in W_P}[E(X_{u})_{w'v}]_P\cdot [\bfE_w(w')]_P\\
&\overset{\text{Lem}.\ref{lem:buw_para}}=\prod_{\al>0,u\al<0}\frac{\theta(\la_{\alv})\theta(\hbar)}{\theta(\la_\alv+\hbar)\theta'(0)}\sum_{w'\in W^P, v\in W_P}
[E(X_{u})_{w'v}]_P\cdot [\bfE_w(w'v)]_P\\
&\overset{\eqref{eq:6}}=\de_{w,u}^{\mathsf{Kr}}. 
\end{align*}
\end{proof}

\end{document}